\documentclass{tglat2e}
\usepackage{amsmath,amssymb,amsfonts,verbatim}

\theoremstyle{plain}
\newtheorem{theorem}{Theorem}
\newtheorem{lemma}[theorem]{Lemma}
\newtheorem{proposition}[theorem]{Proposition}
\newtheorem{corollary}[theorem]{Corollary}
\def\A{\mathcal{E}}
\def\Ab{\bar{\rm A}}
\def\al{\alpha}
\def\Ar{{\rm A}}
\def\At{\widetilde A}

\def\be{\beta}
\def\bt{\ts\,\raise-0.5pt\hbox{\small$\boxtimes$}\,\,}
\def\Bt{\widetilde B}

\def\CC{\mathbb{C}}
\def\com{\ts,\hskip-.5pt}
\def\Cr{{\rm C}}
\def\Ct{\widetilde C}

\def\d{\partial}
\def\de{\delta}
\def\De{\Delta}
\def\Dt{\widetilde D}

\def\End{\operatorname{End}\ts}
\def\ep{\varepsilon}
\def\Ep{E^{\ts\prime}}

\def\ga{\gamma}
\def\ge{\geqslant}
\def\gl{\mathfrak{gl}}

\def\h{\mathfrak{h}}
\def\H{\mathfrak{H}}
\def\Hom{\operatorname{Hom}}

\def\I{{\rm I}}
\def\Ib{\tts\overline{\nns\rm I\nns}\tts}
\def\id{{\rm id}}
\def\io{\iota}

\def\J{{\rm J}}
\def\Jb{\,\overline{\!\rm J\ns}\ts}
\def\Jp{{\rm J}^{\ts\prime}}
\def\Jpb{\,\overline{\!\rm J\ns}\ts^{\ts\prime}}

\def\la{\lambda}
\def\lap{\la^{\ast}}
\def\lcd{\ts,\ldots,}
\def\le{\leqslant}

\def\mup{\mu^{\ast}}

\def\n{\mathfrak{n}}
\def\np{\n^{\ts\prime}}
\def\ns{\hskip-1pt}
\def\nns{\hskip-.5pt}

\def\om{\omega}

\def\op{\oplus}
\def\ot{\otimes}

\def\p{\mathfrak{p}}
\def\P{\mathcal{P}}
\def\Pb{\bar P}
\def\PD{\mathcal{PD}}

\def\q{\mathfrak{q}}
\def\Qb{\bar Q}
\def\qp{\q^{\ts\prime}}

\def\R{{\rm R}}
\def\Rb{\bar{\rm R}}

\def\S{\operatorname{S}}
\def\si{\sigma}
\def\Sym{\mathfrak{S}}

\def\ts{\hskip1pt}
\def\tts{\hskip.5pt}

\def\U{\operatorname{U}}

\def\vp{\varpi}

\def\xic{\check\xi}

\def\Y{\operatorname{Y}}

\def\Z{\operatorname{Z}}
\def\ZZ{{\mathbb Z}} 

\begin{document} 
 
\title[Yangians and Mickelsson Algebras]{Yangians and Mickelsson Algebras I}

\authors{
Sergey Khoroshkin
\address
Institute for Theoretical and\\Experimental Physics\\Moscow 117259, Russia
\email\texttt{khor@itep.ru}
\and 
Maxim Nazarov
\address
Department of Mathematics\\University of York\\York YO10 5DD, England
\email\texttt{mln1@york.ac.uk}
}

\received{September 25, 2005}
\accepted{March 9, 2006}

\maketitle

%------------------------------------------------------------------------------

\renewcommand{\theequation}{\thesection.\arabic{equation}} 
\renewcommand{\thetheorem}{\thesection.\arabic{theorem}}
\makeatletter                                          
\@addtoreset{equation}{section}                        
\makeatother    

%------------------------------------------------------------------------------

\begin{abstract}
We study the composition of the functor
from the category of modules over the Lie algebra $\gl_m$
to the category of modules over the degenerate affine Hecke algebra
of $GL_N$ introduced by I.\,Cherednik,
% T.\,Arakawa, T.\,Suzuki and A.\,Tsuchiya,
with the functor from the latter category to
the category of modules over the Yangian $\Y(\gl_n)$
due to V.\,Drinfeld.
We propose a representation theoretic 
explanation of a link between the intertwining operators on the
tensor products of $\Y(\gl_n)$-modules, and
the \lq\lq\ts extremal cocycle\ts\rq\rq\ 
on the Weyl group of $\gl_m$ defined by D.\,Zhelobenko.
% This link was established by
% by V.\,Tarasov and A.\,Varchenko.
We also establish a connection between 
the composition of the functors, and
the \lq\lq\ts centralizer construction\ts\rq\rq\
of the Yangian $\Y(\gl_n)$ discovered by G.\,Olshanski.
\end{abstract}

%------------------------------------------------------------------------------

\section*{\normalsize 0.\ Introduction}

The central role in this article is played by two well
known constructions. One of these constructions 
is due to V.\,Drinfeld \cite{D2}. Let $\H_N$ be the
degenerate affine Hecke algebra corresponding to the
general linear group $GL_N$ over a non-Archimedean local field. This
is an associative algebra over the complex field $\CC$
which contains the symmetric group ring $\CC\,\Sym_N$ as a subalgebra.
Let $\Y(\gl_n)$ be the Yangian of the general linear Lie algebra $\gl_n\ts$.
This Yangian is a Hopf algebra over the field $\CC$
which contains the universal enveloping algebra $\U(\gl_n)$ as a subalgebra.
In \cite{D2} for any $\H_N$-module $W\ts$, an action of the algebra
$\Y(\gl_n)$ was defined on the vector space
$(W\ot(\CC^{\ts n})^{\ot N}){}^{\ts\Sym_N}$
of the diagonal $\Sym_N$-invariants in the tensor 
product of the vector spaces $W$ and $(\CC^{\ts n})^{\ot N}\ts$.
Details of this construction are reproduced in
Section~1 of the present article. 

The other construction that we use here is due to I.\,Cherednik \cite{C2},
it was also studied by
T.\,Arakawa, T.\,Suzuki and A.\,Tsuchiya \cite{AST}. 
For any module $V$ over the Lie algebra $\gl_m\ts$,
it provides an action of the algebra $\H_N$ on the
tensor product of $\gl_m$-modules $V\ot(\CC^{\ts m})^{\ot N}\ts$. 
This action of $\H_N$
commutes with the diagonal action of $\gl_m$ on the tensor product.
Details of this construction are also reproduced in
Section~1 of the present article. 
By applying the construction from \cite{D2} to
the $\H_N$-module $W=V\ot(\CC^{\ts m})^{\ot N}\ts$,
we get an action of the Yangian $\Y(\gl_n)$ on
$$
(V\ot(\CC^{\ts m})^{\ot N}\ot(\CC^{\ts n})^{\ot N}){}^{\ts\Sym_N}
=
V\ns\ot\ts\S^N(\CC^{\ts m}\ot\CC^{\ts n})
$$
commuting with the action of the Lie algebra $\gl_m\ts$;
see our Proposition \ref{dast}.
By taking the direct sum over $N=0,1,2,\ts\ldots$ of these 
$\Y(\gl_n)$-modules, we turn into an $\Y(\gl_n)$-module
the vector space $V\ns\ot\ts\S\ts(\CC^{\ts m}\ot\CC^{\ts n})\ts$.
It is also a $\gl_m$-module\ts; we denote this bimodule by
$\A_{\ts m}(V)\ts$. The additive group $\CC$ acts on 
the Hopf algebra $\Y(\gl_n)$ by automorphisms.
We denote by $\A_{\ts m}^{\,z}\ts(V)$ the $\Y(\gl_n)$-module obtained
from  $\A_{\ts m}(V)$ via pull-back through the automorphism
of  $\Y(\gl_n)$ corresponding to the element $z\in\CC\ts$.
As a $\gl_{\ts m}\ts$-module $\A_{\ts m}^{\,z}\ts(V)$
coincides with $\A_{\ts m}(V)\ts$.
In this article, we identify the symmetric algebra
$\S\ts(\ts\CC^{\ts m}\ot\CC^{\ts n})$
with the ring $\P\ts(\ts\CC^{\ts m}\ot\CC^{\ts n})$
of polynomial functions on the vector space $\CC^{\ts m}\ot\CC^{\ts n}$.

Now take the Lie algebra $\gl_{m+l}\ts$.
Let $\p$ be the maximal parabolic subalgebra of $\gl_{m+l}$
containing the direct sum of Lie algebras $\gl_m\op\gl_{\ts l}\ts$.
Let $\q$ be the Abelian subalgebra of $\gl_{m+l}$ such that 
$
\gl_{m+l}=\q\op\p\ts.
$
For any $\gl_{\ts l}$-module $U$ let $V\bt U$ be the
$\gl_{m+l}$-module 
parabolically induced from the 
$\gl_m\op\gl_{\ts l}$-module $V\ot U\ts$.
This is a module induced from the subalgebra $\p\ts$.
Consider the space $\A_{\ts m+l}\ts(\ts V\ns\bt U\ts)_{\ts\q}$
of $\q$-coinvariants of the $\gl_{m+l}$-module 
$\A_{\ts m+l}\ts(\ts V\ns\bt U\ts)\ts$.
This space is an $\Y(\gl_n)$-module, which also
inherits the action of the Lie algebra $\gl_m\op\gl_{\ts l}\ts$.
Our Theorem \ref{parind} states that
the bimodule $\A_{\ts m+l}\ts(\ts V\ns\bt U\ts)_{\ts\q}$
over $\Y(\gl_n)$ and $\gl_m\op\gl_{\ts l}$
is equivalent to the tensor product 
$\A_{\ts m}(V)\ot\A_{\ts l}^{\ts m}\ts(U)\ts$.
Here we use the comultiplication on $\Y(\gl_n)\ts$.
The correspondence $V\mapsto\A_{\ts m}(V)$ for $m=1,2,\ts\ldots$
was studied by T.\,Arakawa \cite{A}, but this result seems to be new.
Unlike in \cite{A}, in the proof of Theorem~\ref{parind} we do not use
the representation theory of affine Hecke algebras. 

In Section 3 we propose a representation theoretic
explanation of the correspondence between intertwining operators
on the tensor products of certain $\Y(\gl_n)$-modules, and the
\lq\lq\ts extremal cocycle\ts\rq\rq\ 
on the Weyl group $\Sym_m$ of the reductive Lie algebra
$\gl_m\ts$, defined by D.\,Zhelobenko \cite{Z1}.
This correspondence, discovered
by V.\,Tarasov and A.\,Varchenko \cite{TV2}, was one of the motivations
of our work.
The arguments of \cite{TV2}, inspired by
the results of V.\,Toledano-Laredo \cite{T}, are based on
the classical duality theorem \cite{H}
which asserts that the images of %the algebras
$\U(\gl_m)$ and $\U(\gl_n)$ in the ring
$\PD\ts(\CC^{\ts m}\otimes\CC^{\ts n})$
of differential operators on $\CC^{\ts m}\otimes\CC^{\ts n}$
with polynomial coefficients
are the commutants of each other.
Relevant results were
obtained by Y.\,Smirnov and V.\,Tolstoy \cite{ST}.
Our explanation of the correspondence is based on the theory of
Mickelsson algebras \cite{M1} developed~in~\cite{KO}.
Consider the algebra
\begin{equation}
\label{tenprod}
\U(\gl_m)\ot\PD\ts(\CC^{\ts m}\ot\CC^{\ts n})\ts.
\end{equation}
We have a representation
$\ga:\,\U(\gl_m)\to\PD\ts(\CC^{\ts m}\ot\CC^{\ts n})\ts$.
Taking the composition of the comultiplication map
on $\U(\gl_m)$ with the homomorphism $\id\ot\ga$
we get an embedding of $\U(\gl_m)$ into the algebra \eqref{tenprod}.
Our particular Mickelsson algebra is determined by the pair
formed by the algebra \eqref{tenprod}, and its subalgebra $\U(\gl_m)$
relative to this embedding. %see \cite{K} for more details.
{}From another perspective, connections between the 
representation theory of the Yangian $\Y(\gl_n)$ and 
the theory of Mickelsson algebras have been studied by A.\,Molev \cite{M}.

Our results can be restated in the language of 
dynamical Weyl groups as used by P.\,Etingof and A.\,Varchenko in \cite{EV}.
However, our approach makes more natural the appearance of the Yangian
$\Y(\gl_n)$ in the context of the classical dual pair
$(\gl_m\com\gl_n)$ of reductive Lie algebras.
Moreover, results of the present article
can be extended to other reductive dual pairs \cite{H}.
This will be done in our forthcoming publications.

We complete this article with an observation on
the \lq\lq\ts centralizer construction\ts\rq\rq\ 
of the Yangian $\Y(\gl_n)$ due to G.\,Olshanski \cite{O1}.
For any two irreducible polynomial modules
$V$ and $V^{\ts\prime}$ over the Lie algebra $\gl_m\ts$,
the results of \cite{O1} provide an action
of $\Y(\gl_n)$ on the vector space
$$
\Hom_{\,\gl_m}(\ts V^{\ts\prime}\ts,
V\ot\S\ts(\ts\CC^{\ts m}\ot\CC^{\ts n}\ts))\ts.
$$
Moreover, this action is irreducible. Our Proposition \ref{arol}
states that the same action is inherited from the bimodule 
$\A_{\ts m}(V)=V\ns\ot\ts\S\ts(\CC^{\ts m}\ot\CC^{\ts n})$
over $\Y(\gl_n)$ and $\gl_m\ts$.

We are very grateful to V.\,Tarasov and V.\,Tolstoy for friendly
discussions which led us to the results presented in this article.
The first author was supported by the RFBR grant 05-01-01086, the
grant for Support of Scientific Schools 1999-2003-2, and by
the French National Research Agency grant 
\lq\lq\ts Geometry and Integrability
in Mathematical Physics\ts\rq\rq.
The second author has been supported by the EPSRC grant C511166,
and by the EC grant MRTN-CT2003-505078.
This work was done when both authors stayed at the
Max Planck Institute for Mathematics in Bonn.
We are grateful to the staff of the institute for
their kind help and generous hospitality.

%------------------------------------------------------------------------------

\section*{\normalsize 1.\ Drinfeld functor}
\setcounter{section}{1}
\setcounter{equation}{0}
\setcounter{theorem}{0}

We begin with recalling two well known constructions
from the representation theory of the 
\textit{degenerate affine Hecke algebra\/} $\H_N\ts$, which
corresponds to the general linear group $GL_N$ over a local
non-Archimedean field. This algebra was introduced by V.\,Drinfeld
[D2], see also [L]. By definition, the complex associative algebra $\H_N$ is
generated by the symmetric group algebra $\CC\ts\Sym_N$ and by the pairwise
commuting elements $x_1\lcd x_N$ with the cross relations for $p=1\lcd N-1$
and $q=1\lcd N$
\begin{eqnarray}
\label{cross1}
\si_{p}\,x_q&=&x_q\,\si_{p}\ts,\quad q\neq p\ts,p+1\ts;
\\
\label{cross2}
\si_{p}\,x_p&=&x_{p+1}\,\si_{p}-1\ts.
\end{eqnarray}
Here and in what follows $\si_p\in \Sym_N$ denotes the 
transposition of numbers $p$ and $p+1\ts$. 
More generally, $\si_{pq}\in \Sym_N$ will denote the 
transposition of the numbers $p$ and $q\ts$. The group algebra
$\CC\ts\Sym_N$ can be then regarded as a subalgebra in $\H_N\ts$.
Furhtermore, it follows from the defining relations of $\H_N$
that a homomorphism $\H_N\to\CC\ts\Sym_N\ts$, identical
on the subalgebra $\CC\ts\Sym_N\subset \H_N\ts$,
can be defined by the assignments
\begin{equation}
\label{evalhecke}
x_p\mapsto\si_{1p}+\ldots+\si_{p-1,p}
\quad\text{for}\quad
p=1\lcd N.
\end{equation}

We will also use the elements of the algebra $\H_N\ts$,
$$
y_p=x_p-\si_{1p}-\ldots-\si_{p-1,p}
\quad\text{where}\quad
p=1\lcd N.
$$
Notice that $y_p\mapsto0$ under the homomorphism $\H_N\to\CC\ts\Sym_N\ts$,
defined by
\eqref{evalhecke}. For any permutation $\si\in\Sym_N$, we have
\begin{equation}
\label{yrel}
\si\,y_p\,\si^{-1}=y_{\si(p)}\ts.
\end{equation}
It suffices to verify \eqref{yrel} when
$\si=\si_q$ and $q=1\lcd N-1$. Then (\ref{yrel})
is equivalent to the relations \eqref{cross1},\eqref{cross2}.
The elements $y_1\lcd y_N$ do not commute, but satisfy
the commutation relations
\begin{equation}
\label{ycom}
[y_p\ts,y_q]=\si_{pq}\cdot(y_p-y_q)\ts.
\end{equation}
Let us verify the equality in \eqref{ycom}.
Both sides of \eqref{ycom} are antisymmetric in $p$ and $q\ts$,
so it suffices to consider only the case when $p<q\ts$. Then
\begin{eqnarray}
\nonumber
&&[y_p\ts,x_q]=[x_p-\si_{1p}-\ldots-\si_{p-1,p}\ts,x_q]=0\ts,
\\[4pt]
\nonumber
&&[y_p\ts,y_q\ts]=[y_p\ts,y_q-x_q]=-\ts[y_p\ts,\si_{1q}+\ldots+\si_{q-1,q}]=
-\ts[y_p\ts,\si_{pq}]=\si_{pq}\cdot(y_p-y_q)
\end{eqnarray}
where we used \eqref{yrel}.
Obviously, the algebra $\H_N$ is generated by $\CC\ts\Sym_N$
and the elements $y_1\lcd y_N\ts$.
Together with relations in $\CC\ts\Sym_N$, 
\eqref{yrel} and \eqref{ycom} are defining relations for $\H_N\ts$.
For more details on this presentation of the algebra $\H_N$ see 
%\cite[Section 1.3]{A} and
\cite[Section 5]{N1}.

The first construction we recall here is due to I.\,Cherednik
\cite[Example 2.1]{C2}. It was further studied by
T.\,Arakawa, T.\,Suzuki and A.\,Tsuchiya \cite[Section 5.3]{AST}. 
Let $V$ be any module over the complex general linear Lie algebra $\gl_m\ts$.
Let $E_{ab}\in\gl_m$ with $a,b=1\lcd m$ be the
standard matrix units. We will also regard the matrix units $E_{ab}$
as elements of the algebra $\End(\CC^{\ts m})$, this should not cause any 
confusion. Let us consider the tensor product $V\ot(\CC^{\ts m})^{\ot N}$
of $\gl_m$-modules. Here each of the $N$ tensor factors $\CC^{\ts m}$ is
a copy of the natural $\gl_m$-module. We will use the indices
$1\lcd N$ to label these $N$ tensor factors.
%%% the tensor factor $V$ will be labelled by the index $0$. 
For any $p=1\lcd N$
denote by $E^{\ts(p)}_{ab}$ the operator on the vector space 
$(\CC^{\ts m})^{\ot N}$ acting as
$$
\id^{\ts\ot\ts(p-1)}\ot E_{ab}\ot\id^{\ts\ot\ts(N-p)}\,.
$$

\begin{proposition}
\label{ast}
{\rm\ (i)} 
An action of the algebra $\H_N$ on the vector space 
$V\ot(\CC^{\ts m})^{\ot N}$
can be defined as follows: 
the group $\Sym_N\subset \H_N$ acts naturally by
permutations of the $N$ tensor factors\/ $\CC^{\ts m}$, while any element 
$y_p\in \H_N$ acts as 
\begin{equation}
\label{yact}
\sum_{a,b=1}^m E_{\ts ba}\ot E_{ab}^{\ts(p)}.
\end{equation}
{\rm(ii)} 
This action of $\H_N$ commutes with the
{\rm(}diagonal\,{\rm)} action of\/ 
$\gl_m$ on $V\ot(\CC^{\ts m})^{\ot N}\ts$.
\end{proposition}

\begin{proof}
The relation \eqref{yrel} is obviously satisfied in this
representation of the algebra $\H_N$ by operators on 
$V\ot(\CC^{\ts m})^{\ot N}\ts$.
To verify the relation \eqref{ycom} in this representation, let us
denote by $Y_p$  the operator \eqref{yact}.
For any indices $p,q=1\lcd N$ such that $p\neq q$, then
$$
\begin{aligned}
&[Y_p,Y_q]
=\!\!\sum_{a,b,c,d=1}^m
\big[\,E_{\ts ba}\ot E_{ab}^{\ts(p)},E_{dc}\ot E_{cd}^{\ts(q)}\,\big]
=\!\!\sum_{a,b,c,d=1}^m
\big(\,\de_{ad}\,E_{\ts bc}-\de_{bc}\,E_{da}\big)\ot
E_{ab}^{\ts(p)}\,E_{cd}^{\ts(q)}
\\
&\phantom{[Y_p,Y_q]}
=
\sum_{a,b,c=1}^m E_{\ts bc}\ot E_{ab}^{\ts(p)}\,E_{ca}^{\ts(q)}
-\!\!
\sum_{a,b,d=1}^m E_{da}\ot E_{ab}^{\ts(p)}\,E_{\ts bd}^{\ts(q)}\,.
\end{aligned}
$$
Multiplying on the left
the sum over $a,b,c$ in the last line by the operator
on $V\ot(\CC^{\ts m})^{\ot N}$ corresponding to the transposition 
$\si_{pq}\in\Sym_N$
results in exchanging the first indices $a$ and $c$ in the factors
$E_{ab}^{\ts(p)}$ and $E_{ca}^{\ts(q)}$ in every summand. Since
$$
\sum_{a=1}^m E_{aa}^{\ts(q)}=\id,
$$
the resulting sum equals
$$
\sum_{b,c=1}^m E_{\ts bc}\ot E_{cb}^{\ts(p)}=Y_p\ts.
$$
Similarly, by multiplying on the left the sum over $a,b,d$ by the operator
corresponding to $\si_{pq}\in\Sym_N$ we get the operator $Y_q\ts$. 
This completes the proof of part (i).

The proof of part (ii) consists of a direct verification
that for any $p=1\lcd N$ and $c,d=1\lcd m$
the operator \eqref{yact} commutes with the sum
$$
E_{cd}\ot\id+\id\ot\big(\ts E_{cd}^{\ts(1)}+\ldots+E_{cd}^{\ts(N)}\ts\big)
$$
which describes the action of element $E_{cd}\in\gl_m$
on $V\ot(\CC^{\ts m})^{\ot N}\ts$. Here we omit the verification.
Commutation of the actions of $\Sym_N$ and $\gl_m$ on 
$V\ot(\CC^{\ts m})^{\ot N}$ is obvious.
\qed
\end{proof}

Let us now consider the \textit{triangular decomposition}
of the Lie algebra $\gl_m\ts$,
\begin{equation}
\label{tridec}
\gl_m=\n\op\h\op\np\ts.
\end{equation}
Here $\h$ is the Cartan subalgebra of $\gl_m$ with the basis vectors
$E_{11}\lcd E_{mm}\ts$. Further,
$\n$ and $\np$ are the nilpotent
subalgebras spanned respectively by the elements $E_{\ts ba}$ and
$E_{ab}$ for all $a,b=1\lcd m$ such that $a<b\ts$.
For any $\gl_m$-module $W$, we denote by
$W_{\n}$ the vector space
$W/\ts\n\cdot W$
of the coinvariants of the action of the subalgebra $\n\subset\gl_m$ on $W$.
Note that the Cartan subalgebra $\h\subset\gl_m$ acts on 
the vector space $W_{\n}\ts$.

Now consider the tensor product 
$W=V\ot(\CC^{\ts m})^{\ot N}$ as a left module over the algebra $\H_N\ts$.
The action of $\H_N$ on this module commutes with the
action of the Lie algebra $\gl_m\ts$, 
and hence with the action of the subalgebra $\n\subset\gl_m\ts$.
So the space $(\ts V\ot(\CC^{\ts m})^{\ot N})_{\ts\n}$
of coinvariants of the action of $\n$ 
is a quotient of the $\H_N$-module $V\ot(\CC^{\ts m})^{\ot N}$.
Thus we have a functor from the category of all 
$\gl_m$-modules to the category of bimodules over $\h$ and $\H_N\ts$,
\begin{equation}
\label{zelefun}
V\mapsto(\ts V\ot(\CC^{\ts m})^{\ot N})_{\ts\n}\ts.
\end{equation}
This is a special case of the general construction
due to A.\,Zelevinsky \cite{Z2}.
Restriction of the functor \eqref{zelefun}
to the category $\cal O$
of $\gl_m$-modules \cite{BGG} has been used by
T.\,Arakawa and T.\,Suzuki \cite{AS, S1,S2}
to give algebraic proofs of the Zelevinsky conjecture~\cite{Z1} 
on the multiplicities of composition factors in
the standard modules over the algebra $\H_N\ts$,
and of the Rogawski conjecture~\cite{R} 
on the Jantzen filtration on these modules.
The first of the two conjectures was initially proved
by V.\,Ginzburg \cite{G} using geometric methods.

Now consider the \textit{Yangian}
$\Y(\gl_n)$ of the general linear Lie algebra $\gl_n\ts$.
The Yangian $\Y(\gl_n)$ is a deformation of the universal
enveloping algebra of the polynomial current Lie algebra $\gl_n[u]$
in the class of Hopf algebras, see for instance \cite{D1}.
The unital associative algebra $\Y(\gl_n)$ has a family of generators 
$$
T_{ij}^{\ts(1)},T_{ij}^{\ts(2)},\ts\ldots
\quad\text{where}\quad
i\ts,\ns j=1\lcd n\ts. 
$$
The defining relations for these generators
can be written in terms of the formal series
\begin{equation}
\label{tser}
T_{ij}(u)=
\de_{ij}+T_{ij}^{\ts(1)}u^{-\ns1}+T_{ij}^{\ts(2)}u^{-\ns2}+\,\ldots
\,\in\,\Y(\gl_n)\,[[\ts u^{-1}\ts]]\,.
\end{equation}
Here $u$ is the formal parameter. Let $v$ be another formal parameter.  
Then the defining relations in the associative algebra $\Y(\gl_n)$
can be written as
\begin{equation}
\label{yangrel}
(u-v)\cdot[\ts T_{ij}(u)\ts,T_{kl}(v)\ts]\ts=\;
T_{kj}(u)\ts T_{il}(v)-T_{kj}(v)\ts T_{il}(u)\,,
\end{equation}
where $i\com j\com k\com l=1\lcd n\ts$.
If $n=1$, the algebra $\Y(\gl_n)$ is commutative.
The relations \eqref{yangrel} imply that for any $z\in\CC\,$, the assignments
\begin{equation}
\label{tauz}
\tau_z\ts:\,T_{ij}(u)\ts\mapsto\,T_{ij}(u-z)
\quad\textrm{for}\quad
i\com j=1\lcd n
\end{equation}
define an automorphism $\tau_z$ of the algebra $\Y(\gl_n)\ts$. 
Here each of the formal
power series $T_{ij}(u-z)$ in $(u-z)^{-1}$ should be re-expanded 
in $u^{-1}$, and the assignment \eqref{tauz} is a correspondence
between the respective coefficients of series in $u^{-1}$.

Now let $E_{ij}\in\gl_n$ with $i,j=1\lcd n$ be the
standard matrix units. We will also regard the matrix units $E_{ij}$
as elements of the algebra $\End(\CC^{\ts n})$, this should not cause any 
confusion. The Yangian $\Y(\gl_n)$ contains 
the universal enveloping algebra $\U(\gl_n)$ as a subalgebra\ts;
the embedding $\U(\gl_n)\to\Y(\gl_n)$ can be defined by the assignments
$$
%\label{uemb}
E_{ij}\mapsto T_{ij}^{\ts(1)}
\quad\text{for}\quad
i\ts,\ns j=1\lcd n\ts. 
$$
Moreover, there is a homomorphism $\pi_n:\Y(\gl_n)\to\U(\gl_n)$ identical
on the subalgebra $\U(\gl_n)\subset\Y(\gl_n)\ts$, 
it can be defined by the assignments
\begin{equation}
\label{pin}
\pi_n\ts:\,T_{ij}^{\ts(2)},T_{ij}^{\ts(3)},\ts\ldots
\,\mapsto\,0
\quad\text{for}\quad
i\ts,\ns j=1\lcd n\ts. 
\end{equation}
For further details on the definition of 
the algebra $\Y(\gl_n)$ see \cite[Chapter 1]{MNO}.

The second construction we recall here is due to V.\,Drinfeld \cite{D2},
this construction has originally motivated his 
definition of the degenerate affine Hecke
algebra $\H_N\ts$. 
For $p=1\lcd N$
denote by $E^{\ts(p)}_{ij}$ the operator on the vector space 
$(\CC^{\ts n})^{\ot N}$ acting as
$$
\id^{\ts\ot\ts(p-1)}\ot E_{ij}\ot\id^{\ts\ot\ts(N-p)}\,.
$$
The group $\Sym_N$ acts on the tensor
product $(\CC^{\ts n})^{\ot N}$ from the left
by permutations of the $N$ tensor factors.
Let $W$ be any $\H_N$-module.
The group $\Sym_N$ also
acts from the left on $W$, 
via the embedding $\CC\ts\Sym_N\to \H_N$.
Consider the subspace
\begin{equation}
\label{sinv}
(W\ot(\CC^{\ts n})^{\ot N}){}^{\ts\Sym_N}\subset \,W\ot(\CC^{\ts n})^{\ot N}
\end{equation}
of invariants with respect to the diagonal action of $\Sym_N\ts$.
In the next proposition we use the convention that
$(\ts y_p)^{\ts0}=1$, the identity element of the algebra $\CC\ts\Sym_N\ts$.

\begin{proposition}
\label{d}
One can define
an action of the algebra\/ $\Y(\gl_n)$ on the vector space
$(W\ot(\CC^{\ts n})^{\ot N}){}^{\ts\Sym_N}$ so that 
for any $s=0,1,2,\ts\ldots$ the generator
$T_{ij}^{\ts(s+1)}$ acts as
\begin{equation}
\label{tact}
\sum_{p=1}^N\,(-\ts y_p)^s\ot E^{\ts(p)}_{ij}\ts.
\end{equation}
\end{proposition}

\begin{proof}
As an operator on the vector space $W\ot(\CC^{\ts n})^{\ot N}$,
\eqref{tact}
commutes with the diagonal action of the group $\Sym_N\ts$, due to the
relations \eqref{yrel} for the generators $y_1\lcd y_N$ of $\H_N\ts$.
So the restriction of the operator \eqref{tact} to the
subspace \eqref{sinv} is well defined. 

By substituting the sum
\eqref{tact} for every $T_{ij}^{\ts(s+1)}$
in \eqref{tser}, we get the series in $u^{-1}$
with the coefficients in the algebra $\H_N\ot\End((\CC^{\ts n})^{\ot N})$, 
$$
\de_{ij}\ot\id\,+
\sum_{s=0}^\infty\,
\sum_{p=1}^N\,
(-\ts y_p)^s\,u^{-s-1}\ot
E^{\ts(p)}_{ij}
=\ts
\de_{ij}\ot\id\,+
\sum_{p=1}^N\,
(u+y_p)^{-1}\ot
E^{\ts(p)}_{ij}
\ts.
$$
Making the respective substitutions for $T_{ij}(u)$ and $T_{kl}(v)$
at the left hand side of the defining relations \eqref{yangrel},
and then cancelling the commutators with $\de_{ij}\ot\id$ and 
$\de_{kl}\ot\id$, we obtain the sum
$$
(u-v)\sum_{p,q=1}^N\,
\big(\ts
(u+y_p)^{-1}\,(v+y_q)^{-1}\ot
E^{\ts(p)}_{ij}\,E^{\ts(q)}_{kl}
-
(v+y_q)^{-1}\,(u+y_p)^{-1}\ot
E^{\ts(q)}_{kl}\,E^{\ts(p)}_{ij}
\ts\big)
\vspace{-16pt}
$$
\begin{eqnarray}
\label{sum11}
=\ &(u-v)&\,\sum_{p=1}^N\ 
(u+y_p)^{-1}\,(v+y_p)^{-1}\ot
[\ts E^{\ts(p)}_{ij},E^{\ts(p)}_{kl}\ts]\,+
\\
\label{sum12}
\,&(u-v)&\,\sum\limits_{\substack{p,q=1\\p\neq q}}^N\,
[\,(u+y_p)^{-1},(v+y_q)^{-1}\ts]\ot
E^{\ts(p)}_{ij}\,E^{\ts(q)}_{kl}\ts.
\end{eqnarray}
Making the substitutions at the right hand side of \eqref{yangrel},
and cancelling the two tensor products $\de_{kj}\,\de_{il}\ot\id$ in
the resulting difference, we get
\begin{equation}
\label{sum13}
\sum_{p=1}^N\,
\big(\ts
(v+y_p)^{-1}-(u+y_p)^{-1}\ts\big)\ot
\big(\,\de_{kj}\,E^{\ts(p)}_{il}-\de_{il}\,E^{\ts(p)}_{kj}\ts\big)\,\ts+
\end{equation}
\vspace{-8pt}
\begin{equation}
\label{sum14}
\sum_{p,q=1}^N\,
\big(\ts
(u+y_p)^{-1}\,(v+y_q)^{-1}
-
(v+y_p)^{-1}\,(u+y_q)^{-1}
\ts\big)
\ot
E^{\ts(p)}_{kj}\,E^{\ts(q)}_{il}\ts.
\end{equation}
The sums \eqref{sum11} and \eqref{sum13} are equal to each other.
In the sum \eqref{sum14}, the summands with $p=q$ vanish. In every summand
of \eqref{sum14} with $p\neq q$, the factors 
$E^{\ts(p)}_{kj}$ and $E^{\ts(q)}_{il}$ commute. Hence,
by exchanging the indices $p$ and $q\ts$, the sum \eqref{sum14} equals
$$
\sum\limits_{\substack{p,q=1\\p\neq q}}^N\,
\big(\ts
(u+y_q)^{-1}\,(v+y_p)^{-1}
-
(v+y_q)^{-1}\,(u+y_p)^{-1}
\ts\big)
\ot
E^{\ts(p)}_{il}\,E^{\ts(q)}_{kj}\ts.
$$
The action of the latter sum on the subspace \eqref{sinv}
%of $\Sym_N$-invariants 
coincides with the action of the sum
\begin{equation}
\label{sum15}
\sum\limits_{\substack{p,q=1\\p\neq q}}^N\,
\big(\ts
(u+y_q)^{-1}\,(v+y_p)^{-1}
-
(v+y_q)^{-1}\,(u+y_p)^{-1}
\ts\big)\cdot\si_{pq}\ts
\ot
E^{\ts(p)}_{ij}\,E^{\ts(q)}_{kl}\ts.
\end{equation}
The sum \eqref{sum12} is equal to \eqref{sum15}, because for $p\neq q$ 
we have the relation
$$
(u-v)\cdot[\,(u+y_p)^{-1},(v+y_q)^{-1}\ts]=
$$
$$
\big(\ts
(u+y_q)^{-1}\,(v+y_p)^{-1}
-
(v+y_q)^{-1}\,(u+y_p)^{-1}
\ts\big)\cdot\si_{pq}\ts.
$$
To verify this relation, let us multiply its sides by
$(u+y_p)\ts(v+y_q)$ on the left, and by $(v+y_q)\ts(u+y_p)$
on the right. Using the equality
$$
\si_{pq}\cdot(v+y_q)\ts(u+y_p)=(v+y_p)\ts(u+y_q)\cdot\si_{pq}\ts,
$$
then we get the relation
$$
(u-v)\cdot[\ts u+y_p,v+y_q\ts]=
$$
$$
\big(\ts
(u+y_p)\,(v+y_q)
-
(v+y_p)\,(u+y_q)
\ts\big)\cdot\si_{pq}\ts.
$$
But the last relation holds true due to \eqref{yrel} and \eqref{ycom}.
\qed
\end{proof}

\noindent\textit{Remark.}
When $s=0\ts$, the sum \eqref{tact} describes the action
of the element $E_{ij}\in\gl_n$ on the tensor product space 
$W\ot(\CC^{\ts n})^{\ot N}$, and hence on its subspace \eqref{sinv}.
Here each of the $N$ tensor factors $\CC^{\ts n}$ is regarded as a copy
of the natural $\gl_n$-module, and the action of $\gl_n$ on $W$ is trivial.
Hence the action of the Yangian $\Y(\gl_n)$ on the subspace
\eqref{sinv} as defined in Proposition \ref{d},
is compatible with the embedding $\U(\gl_n)\to\Y(\gl_n)\ts$.
\qed

\medskip\smallskip
Thus we obtain a functor from the category of all 
$\H_N$-modules to the category of $\Y(\gl_n)$-modules
\begin{equation}
\label{drinfun}
W\mapsto(W\ot(\CC^{\ts n})^{\ot N}){}^{\ts \Sym_N}\ts.
\end{equation}
This is the \textit{Drinfeld functor} for the 
Yangian $\Y(\gl_n)\ts$.
Let us now apply this functor to the $\H_N$-module
$W=V\ot(\CC^{\ts m})^{\ot N}$ where $V$ is an arbitrary $\gl_m$-module;
see Proposition~\ref{ast}. The vector space of the resulting
$\Y(\gl_n)$-module is
$$
(V\ot(\CC^{\ts m})^{\ot N}\ot(\CC^{\ts n})^{\ot N}){}^{\ts\Sym_N}
=
V\ot((\CC^{\ts m}\ot\CC^{\ts n})^{\ot N}){}^{\ts\Sym_N}
$$
where the group $\Sym_N$ acts by permutations of
the $N$ tensor factors $\CC^{\ts m}\ot\CC^{\ts n}\ts$. Hence the
resulting vector space is 
\begin{equation}
\label{yangmod}
V\ns\ot\ts\S^N(\CC^{\ts m}\ot\CC^{\ts n})
\end{equation}
where we take the $N$-th symmetric power of the vector space 
$\CC^{\ts m}\ot\CC^{\ts n}\ts$. Note that the Lie algebra $\gl_m$ also acts
on \eqref{yangmod} as the tensor product of two $\gl_m$-modules. 

We can identify the vector space $\CC^{\ts m}\ot\CC^{\ts n}$
with its dual vector space, 
so that the standard basis vectors of $\CC^{\ts m}\ot\CC^{\ts n}$
are identified with the corresponding coordinate functions
$x_{ai}$ 
where
$a=1\lcd m$
and
$i=1\lcd n\ts$.
The symmetric algebra $\S\ts(\CC^{\ts m}\ot\CC^{\ts n})$
is then identified with the ring $\P\ts(\CC^{\ts m}\ot\CC^{\ts n})$
of polynomial functions on $\CC^{\ts m}\ot\CC^{\ts n}\ts$.
The ring of differential operators on $\P\ts(\CC^{\ts m}\ot\CC^{\ts n})$
with polynomial coefficients will be denoted by 
$\PD\ts(\CC^{\ts m}\ot\CC^{\ts n})\ts$.
Let $\d_{ai}$ be the partial derivation on $\P\ts(\CC^{\ts m}\ot\CC^{\ts n})$
corresponding to the variable $x_{ai}\ts$.
We can now describe the action of $\Y(\gl_n)$ on
the vector space \eqref{yangmod}.
%cf.\ \cite[Section 3]{A}.
%by combining Propositions \ref{ast} and \ref{d}.

\begin{proposition}
\label{dast}
{\rm\,(i)} 
For any $s=0,1,2,\ts\ldots$ the generator
$T_{ij}^{\ts(s+1)}$
acts on the $\Y(\gl_n)$-module \eqref{yangmod}
as the element of the tensor product\/ 
$\U(\gl_m)\ot\PD\ts(\CC^{\ts m}\ot\CC^{\ts n})\ts$,
\begin{equation}
\label{combact}
\sum_{c_0,c_1,\ldots,c_s=1}^m
(-1)^s\,E_{c_1c_0}\,E_{c_2c_1}\ldots\ts E_{c_sc_{s-1}}
\ot x_{c_0i}\,\d_{c_sj}\ts.
\end{equation}
In the case $s=0$, the first tensor factor in the summand in\/
\eqref{combact} is understood as $1$.
{\rm\ (ii)} 
The action of $\Y(\gl_n)$ on \eqref{yangmod}
commutes with the {\rm(}diagonal\,{\rm)} action of\/ $\gl_m\ts$.
\end{proposition}

\begin{proof}
First let us consider the action of the sum \eqref{tact} on the
vector space $W\ot(\CC^{\ts n})^{\ot N}$ where 
$W=V\ot(\CC^{\ts m})^{\ot N}$\ts.
By substituting the sum \eqref{yact} for $y_p$ in \eqref{tact},
we then get the sum
\begin{equation}
\label{subst}
\sum_{p=1}^N\,
\Big(-\!
\sum_{a,b=1}^m E_{\ts ba}\ot E_{ab}^{\ts(p)}
\,\Big){\!\!\ns\phantom{\big)}}^s\ot E^{\ts(p)}_{ij}
\end{equation}
acting on the vector space 
$V\ot(\CC^{\ts m})^{\ot N}\ot(\CC^{\ts n})^{\ot N}$.
Using the relations
$$
E_{ab}^{\ts(p)}\,E_{cd}^{\ts(p)}=\de_{bc}\,E_{ad}^{\ts(p)}
\quad\text{for}\quad
p=1\lcd N
$$
the sum \eqref{subst} can be rewritten as %the sum over $p=1\lcd N$ of
$$
\sum_{p=1}^N\ 
\sum_{c_0,c_1\ldots,c_s=1}^m\ 
(-1)^s\,E_{c_1c_0}
\ldots\ts 
E_{c_sc_{s-1}}\ot
E_{c_0c_s}^{\ts(p)}\ot
E^{\ts(p)}_{ij}\ts.
$$
To prove the part (i) of the proposition, 
it remains to observe that after identifying the subspace
$$
((\CC^{\ts m})^{\ot N}\ot(\CC^{\ts n})^{\ot N}){}^{\ts\Sym_N}
\subset
(\CC^{\ts m})^{\ot N}\ot(\CC^{\ts n})^{\ot N}
$$
with the space $\P^N(\CC^{\ts m}\ot\CC^{\ts n})$
of polynomial functions on $\CC^{\ts m}\ot\CC^{\ts n}$
of degree $N$, the operator
$$
\sum_{p=1}^N\,\,
E_{c_0c_s}^{\ts(p)}\ot
E^{\ts(p)}_{ij}
$$
on this subspace
can be identified with the operator $x_{c_0i}\,\d_{c_sj}$
on the space $\P^{\ts N}(\CC^{\ts m}\ot\CC^{\ts n})\ts$.
The part (ii) of Proposition \ref{dast}
follows from the respective part of Proposition \ref{ast}.
\qed
\end{proof}

\noindent\textit{Remark.}
By definition, the basis element $E_{ab}\in\gl_m$ acts on the
vector space \eqref{yangmod} as 
\begin{equation}
\label{eabact}
E_{ab}\ot1+
\sum_{k=1}^n\,
1\ot x_{ak}\,\d_{\ts bk}\ts.
\end{equation}
One can easily verify by straightforward calculation,
that the elements \eqref{combact} and \eqref{eabact}
of the algebra $\U(\gl_m)\ot\PD\ts(\CC^{\ts m}\ot\CC^{\ts n})$
commute with each other. Moreover, using the First Fundamental Theorem
of invariants for the general linear group $GL_m\ts$,
one can show that the commutant
in the algebra $\U(\gl_m)\ot\PD\ts(\CC^{\ts m}\ot\CC^{\ts n})$
of all elements \eqref{eabact} with $a,b=1\lcd m$ 
is generated by the subalgebra $\Z(\gl_m)\ot1$ and
all elements of the form \eqref{combact};
cf.\ \cite[Section 2.1]{O2}.
Here $\Z(\gl_m)$ denotes
the centre of the universal enveloping algebra $\U(\gl_m)$.
This extends the classical theorem
\cite[Section 2.3]{H} stating that the two families of operators
on the vector space $\P\ts(\CC^{\ts m}\ot\CC^{\ts n})\ts$,
\begin{equation}
\label{glmact}
\sum_{k=1}^n\,
x_{ak}\,\d_{\ts bk}
\quad\text{where}\quad
a,b=1\lcd m
\end{equation}
and
\begin{equation}
\label{glnact}
\sum_{c=1}^m\,
x_{ci}\,\d_{cj}
\quad\text{where}\quad
i,j=1\lcd n
\end{equation}
generate their mutual commutants in the algebra 
$\PD\ts(\CC^{\ts m}\ot\CC^{\ts n})$.
Here the operators \eqref{glmact} and \eqref{glnact}
describe the
actions on $\P\ts(\CC^{\ts m}\ot\CC^{\ts n})$ of the basis elements 
$E_{ab}\in\gl_m$ and $E_{ij}\in\gl_n$ respectively.
\qed

\medskip\smallskip
We finish this section with an observation
on matrices with entries from the universal enveloping
algebra $\U(\gl_m)\ts$.
Let $E$ be the $m\times m$ matrix whose $ab\ts$-entry is 
the generator $E_{ab}\in\gl_m\ts$. Let $\Ep$ be the
transposed matrix.
Take the matrix inverse to $u+\Ep\ts$.
Here the summand $u$ stands for the scalar $m\times m$ matrix
with diagonal entry $u\ts$, and the inverse is a formal power
series in $u^{-1}$ with matrix coefficients. 
Denote by $X_{ab}(u)$ the $ab\ts$-entry of inverse matrix. Then
$$
X_{ab}(u)=
\ts\de_{ab}\,u^{-1}\ns-
E_{\ts ba}\ts u^{-2}\,+
$$
\begin{equation}
\label{xabu}
\sum_{s=1}^\infty\,
\sum_{c_1,\ldots,c_s=1}^m
(-1)^{s+1}\,
E_{c_1a}\ts E_{c_2c_1}\ldots\ts E_{c_sc_{s-1}}\ts E_{\ts bc_s}
u^{-s-2}\ts.
\end{equation}
The assignment of the element \eqref{combact}
to any coefficient $T_{ij}^{\ts(s+1)}$
of the series \eqref{tser} can be now written as
$$
T_{ij}(u)\mapsto\de_{ij}\ts+\sum_{a,b=1}^m
X_{ab}(u)\ot x_{ai}\,\d_{\ts bj}\ts.
$$

%------------------------------------------------------------------------------

\section*{\normalsize 2.\ Parabolic induction}
\setcounter{section}{2}
\setcounter{equation}{0}
\setcounter{theorem}{0}

The Yangian $\Y(\gl_n)$ is a Hopf algebra over the field $\CC\ts$.
Using the series \eqref{tser},
the comultiplication $\De:\Y(\gl_n)\to\Y(\gl_n)\ot\Y(\gl_n)$ is defined by
the assignment
\begin{equation}\label{1.33}
\De:T_{ij}(u)\ts\mapsto\ts\sum_{k=1}^n\ T_{ik}(u)\ot T_{kj}(u)\,;
\end{equation}
the tensor product at the right hand side of the assignment (\ref{1.33})
is taken over the subalgebra 
$\CC[[u^{-1}]]\subset\Y(\gl_n)\,[[u^{-1}]]\ts$.
When taking tensor products of modules over $\Y(\gl_n)$, 
we use the comultiplication \eqref{1.33}.
The counit homomorphism 
$\ep:\Y(\gl_n)\to\CC$ is defined by 
$$
\ep:\,T_{ij}(u)\ts\mapsto\ts\de_{ij}\cdot1\ts.
$$
The antipode ${\rm S}$ on $\Y(\gl_n)$ is defined by using the 
$n\times n$ matrix $T(u)$ whose $ij$-entry is the series $T_{ij}(u)\ts$.
This matrix is invertible as formal power series in $u^{-1}$ with
matrix coefficients, because the leading term of this series is
the identity $n\times n$ matrix. Then
the involutive anti-automorphism ${\rm S}$ of $\Y(\gl_n)$
is defined by the assignment
$$
{\rm S}\ts:\ts T(u)\mapsto T(u)^{-1}.
$$
This assignment means that by applying ${\rm S}$ to the coefficients of
the series $T_{ij}(u)$, we obtain the series which is the
$ij$-entry of the inverse matrix $T(u)^{-1}\ts$. 
We also use the involutive automorphism $\om_n$
of the algebra $\Y(\gl_n)$ defined by a similar assignment,
\begin{equation}
\label{1.51}
\om_n:\ts T(u)\mapsto T(-u)^{-1}.
\end{equation}
For more details on the Hopf algebra
structure on $\Y(\gl_n)$ see again \cite[Chapter 1]{MNO}.

Let us now consider the infinite direct sum of bimodules over
$\gl_m$ and $\Y(\gl_n)\ts$,
$$ 
\mathop{\op}\limits_{N=0}^\infty\,V\ns\ot\ts\S^N(\CC^{\ts m}\ot\CC^{\ts n})
=V\ns\ot\ts\S\,(\CC^{\ts m}\ot\CC^{\ts n})\ts.
$$
Let us denote this bimodule by $\A_{\ts m}(V)\ts$, so that $\A_{\ts m}$
is a functor from the category of all $\gl_m$-modules to the
category of bimodules over $\gl_m$ and $\Y(\gl_n)\ts$. 
By identifying the symmetric algebra
$\S\,(\CC^{\ts m}\ot\CC^{\ts n})$ with the ring 
$\P\ts(\CC^{\ts m}\ot\CC^{\ts n})$, 
the action of the generator $T_{ij}^{\ts(s+1)}$ of $\Y(\gl_n)$
on $\A_{\ts m}(V)$ is described by the formula \eqref{combact}.

For any positive integer $l$ let $U$ be a
module over the Lie algebra $\gl_{\ts l}\ts$. Then $\A_{\ts l}(U)$ 
is another $\Y(\gl_n)$-module. For any $z\in\CC$ denote by
$\A_{\ts l}^{\ts z}\ts(U)$ the $\Y(\gl_n)$-module obtained
from  $\A_{\ts l}(U)$ via pull-back through the automorphism
$\tau_z$ of $\Y(\gl_n)\ts$, defined by \eqref{tauz}. 
As a $\gl_{\ts l}\ts$-module $\A_{\ts l}^{\ts z}\ts(U)$
coincides with $\A_{\ts l}(U)\ts$.

The decomposition $\CC^{\ts m+l}=\CC^{\ts m}\op\CC^{\ts l}$ 
determines an embedding
of the direct sum $\gl_m\op\gl_{\ts l}$ of Lie algebras into $\gl_{m+l}\ts$.
As a subalgebra of $\gl_{m+l}\ts$,
the direct summand $\gl_m$ is spanned by the matrix units 
$E_{ab}\in\gl_{m+l}$ where
$a,b=1\lcd m\ts$. The direct summand $\gl_{\ts l}$ is spanned by
the matrix units $E_{ab}$ where $a,b=m+1\lcd m+l\ts$.
Let $\q$ and $\qp$ be 
the Abelian subalgebras of $\gl_{m+l}$
spanned respectively by matrix units $E_{\ts ba}$ and 
$E_{ab}$ for all $a=1\lcd m$ and $b=m+1\lcd m+l\ts$.
Put $\p=\gl_m\op\gl_{\ts l}\op\qp\ts$.
Then $\p$ is a maximal parabolic subalgebra of the reductive 
Lie algebra $\gl_{m+l}\ts$, and moreover 
$
\gl_{m+l}=\q\op\p\ts.
$
Denote by $V\ns\bt U$ the $\gl_{m+l}$-module \textit{parabolically induced\/}
from the $\gl_m\op\gl_{\ts l\ts}$-module $V\ot U$.
To define $V\ns\bt U$, one first extends the action of the Lie algebra
$\gl_m\op\gl_{\ts l\ts}$ on $V\ot U$ to the Lie algebra $\p\ts$, so that
any element of the subalgebra $\qp\subset\p$ acts on $V\ot U$ as zero.
By definition, $V\ns\bt U$ is the $\gl_{m+l}$-module induced from the 
$\p$-module $V\ot U$.  

Now consider the bimodule $\A_{\ts m+l}\ts(\ts V\ns\bt U\ts)$
over $\gl_{m+l}$ and $\Y(\gl_n)\ts$.
Here the action of $\Y(\gl_n)$ commutes with the
action of the Lie algebra $\gl_{m+l}\ts$, 
and hence with the action of the subalgebra $\q\subset\gl_{m+l}\ts$.
Therefore the vector space 
$\A_{\ts m+l}\ts(\ts V\ns\bt U\ts)_{\ts\q}$
of coinvariants of the action of the subalgebra $\q$ 
is a quotient of the $\Y(\gl_n)$-module 
$\A_{\ts m+l}\ts(\ts V\ns\bt U\ts)\ts$.
Note that the subalgebra 
$\gl_m\op\gl_{\ts l}\subset\gl_{m+l}$ also acts on this quotient space.

\begin{theorem}
\label{parind}
The bimodule $\A_{\ts m+l}\ts(\ts V\ns\bt U\ts)_{\ts\q}$
over the Yangian $\Y(\gl_n)$ and
the direct sum $\gl_m\op\gl_{\ts l}\ts$, 
is equivalent to the tensor product 
$\A_{\ts m}(V)\ot\A_{\ts l}^{\ts m}\ts(U)\ts$.
\end{theorem}

Our proof of the theorem is based on two simple lemmas.
The first of these lemmas 
applies to matrices over arbitrary unital ring.
Take a matrix of size $(m+l)\times(m+l)$ over such a ring,  
and write it as the block matrix
\begin{equation}
\label{blockmat}
\begin{bmatrix}\,A\,&B\,\\\,C\,&D\,\end{bmatrix}
\end{equation}
where the blocks $A,B,C,D$ are matrices of sizes
$m\times m$, $m\times l$, $l\times m$, $l\times l$ respectively. 
The following fact is well known, see for instance \cite[Lemma 3.2]{B}.

\begin{lemma}
\label{lemma1}
Suppose the matrix \eqref{blockmat} is invertible.
Suppose the matrices $A$ and $D$ are also invertible. Then the matrices
$A-B\ts D^{-1}\ts C$ and $D-C\ts A^{-1}B$ are invertible too, and
$$
\begin{bmatrix}A&B\\C&D\end{bmatrix}^{-1}\!=\ \  
\begin{bmatrix}
(A-B\ts D^{-1}\ts C)^{-1}
&
\,-A^{-1}B\ts(\ts D-C\ts A^{-1}B)^{-1}\,
\\
\,-D^{-1}C\ts(A-B\ts D^{-1}\ts C)^{-1}\,
&
(\ts D-C\ts A^{-1}B)^{-1}
\end{bmatrix}
\,.
$$
\end{lemma}

Consider again
the $m\times m$ matrix $E$ whose $ab\ts$-entry is 
the generator $E_{ab}\in\gl_m\ts$. The $ab\ts$-entry $X_{ab}(u)$
of the matrix
inverse to $u+\Ep$ is given by the equality \eqref{xabu}.
Denote by $Z(u)$ the trace of the inverse matrix, so that
\begin{equation}
\label{zu}
Z(u)\,=\,\sum_{c=1}^m\,X_{cc}(u)\ts.
\end{equation}
Then $Z(u)$ is a formal power series in $u^{-1}$
with the coefficients from the algebra  $\U(\gl_m)\ts$.
Note that the leading term of this series is $m\ts u^{-1}$.
Let us now regard the coefficents of the series \eqref{xabu} and
\eqref{zu} as elements of the algebra
$\U(\gl_{m+l})\ts$, using the standard embedding of 
the Lie algebra $\gl_m$ to $\gl_{m+l}\ts$.

\begin{lemma}
\label{lemma2}
For any\/ $a=1\lcd m$ and\/ $d=1\lcd l$ 
we have an equality of the series
with coefficients in the algebra\/ $\U(\gl_{m+l})\ts,$
\begin{equation}
\label{xe}
\sum_{b=1}^m\,E_{\ts m+d,b}\,X_{ab}(u)
\,=\,
\sum_{b=1}^m\,X_{ab}(u)\,E_{\ts m+d,b}\,(1-Z(u))\ts.
\end{equation}
\end{lemma}

\begin{proof}
For any indices $b\com c\com d\com e=1\lcd m$ we have the equality
$$
(\ts\de_{ec}\,u+E_{ec}\ts)\,E_{\ts m+d,b}
\,=\,
E_{\ts m+d,b}\,(\ts\de_{ec}\,u+E_{ec}\ts)-\de_{eb}\,E_{\ts m+d,c}\,.
$$
Multiplying both sides of this equality by $X_{ac}(u)$ on the left
and taking the sums over $c=1\lcd m$ we obtain the equality
$$
\de_{ae}\,E_{\ts m+d,b}
\,=\,\sum_{c=1}^m\,
X_{ac}(u)\,E_{\ts m+d,b}\,(\ts\de_{ec}\,u+E_{ec}\ts)
\,-\,\sum_{c=1}^m\,
X_{ac}(u)\,\de_{eb}\,E_{\ts m+d,c}\,.
$$
Multiplying both sides of the latter equality by $X_{eb}(u)$ on the right
and taking the sums over $e=1\lcd m$ we get the equality
$$
E_{\ts m+d,b}\,X_{ab}(u)
\,=\,
X_{ab}(u)\,E_{\ts m+d,b}
\,-\,\sum_{c=1}^m
X_{ac}(u)\,E_{\ts m+d,c}\,X_{bb}(u)\,.
$$
Taking here the sums over $b=1\lcd m$ and using the definition
\eqref{zu} we get \eqref{xe}.
\qed
\end{proof}

\noindent\textit{Proof of Theorem \ref{parind}.}
The %underlying 
vector space of the $\gl_{m+l}\ts$-module $V\bt U$
can be identified with the tensor product $\U(\q)\ot V\ns\ot U$
so that the Lie subalgebra $\q\subset\gl_{m+l}$ acts via
left multiplication on the first tensor factor.
Note that the corresponding
action of the commutative algebra $\U(\q)$ is free.
The tensor product $V\ns\ot U$ is then identified with the subspace
\begin{equation}
\label{onesub}
1\ot V\ns\ot U\subset \U(\q)\ot V\ns\ot U\ts.
\end{equation}
On this subspace,
any element of the subalgebra $\qp\subset\gl_{m+l}$ acts as zero, 
while the two direct summands
of subalgebra $\gl_m\op\gl_{\ts l}\subset\gl_{m+l}$
act non-trivially only on the tensor factors $V$ and $U$ respectively.
All this determines the action of Lie algebra $\gl_{m+l}$ on
$\U(\q)\ot V\ns\ot U$. Now consider $\A_{\ts m+l}\ts(\ts V\ns\bt U\ts)$
as a $\gl_{m+l}\ts$-module, we will denote it by $W$ for short.
Then $W$ is the tensor product of two $\gl_{m+l}\ts$-modules,
$$
W=(V\bt U)\ot\P\ts(\CC^{\ts m+l}\ot\CC^{\ts n})=
\U(\q)\ot V\ns\ot U\ot\P\ts(\CC^{\ts m+l}\ot\CC^{\ts n})\ts.
$$

The vector spaces of the two $\Y(\gl_n)$-modules 
$\A_{\ts m}(V)$ and $\A_{\ts l}^{\ts m}\ts(U)$ are respectively
$$
V\ns\ot\ts\P\ts(\CC^{\ts m}\ot\CC^{\ts n})
\ \quad\text{and}\ \quad
U\ns\ot\ts\P\ts(\CC^{\ts l}\ot\CC^{\ts n})\ts.
$$
Identify the tensor product of these two vector spaces with
\begin{equation}
\label{vuprod}
V\ns\ot U\ot
\P\ts(\CC^{\ts m}\ot\CC^{\ts n})\ot\P\,(\CC^{\ts l}\ot\CC^{\ts n})=
V\ns\ot U\ot
\P\ts(\CC^{\ts m+l}\ot\CC^{\ts n})
\end{equation}
where we use the standard direct sum decomposition 
$$\CC^{\ts m+l}\ot\CC^{\ts n}=
\CC^{\ts m}\ot\CC^{\ts n}\op\ts\CC^{\ts l}\ot\CC^{\ts n}\ts.
$$
Regard the tensor product $V\ns\ot U$ in \eqref{vuprod} as a module over
the subalgebra $\gl_m\op\gl_{\ts l}\subset\gl_{m+l}$\ts. This subalgebra
also acts on $\P\ts(\CC^{\ts m+l}\ot\CC^{\ts n})$ naturally.
Define a linear map
$$
\chi:\,V\ns\ot U\ot\P\ts(\CC^{\ts m+l}\ot\CC^{\ts n})\ts\to W\ts/\,\q\cdot W
$$
by the assignment
$$
\chi:\,y\ot x\ot f\,\mapsto\,1\ot y\ot x\ot f\,+\,\q\cdot W
$$
for any $y\in V$, $x\in U$ and $f\in\P\ts(\CC^{\ts m+l}\ot\CC^{\ts n})\ts$. 
The operator $\chi$ evidently intertwines the actions of the Lie
algebra $\gl_m\op\gl_{\ts l}\ts$. 

%This operator is surjective, because
%the subspace \eqref{onesub} is $\U(\q)$-cyclic.

Let us demonstrate that the operator $\chi$ is bijective.
Firstly consider the action of the Lie subalgebra $\q\subset\gl_{m+l}$
on the vector space
$$
\P\ts(\CC^{\ts m+l})\ts=\ts
\P\ts(\CC^{\ts m})\ot\P\ts(\CC^{\ts l})\ts.
$$
This vector space admits an ascending filtration by the subspaces
$$
\mathop{\op}\limits_{N=0}^{\ts K}  
\P^{\ts N}(\CC^{\ts m})\ot\P\ts(\CC^{\ts l})
\ \quad\textrm{where}\ \quad
K=0,1,2,\ts\ldots\,.
$$
Here $\P^{\ts N}(\CC^{\ts m})$ is
the space of polynomial functions on $\CC^{\ts m}$ of degree $N\ts$.
The action of the Lie algebra $\q$ on $\P\ts(\CC^{\ts m+l})$
preserves each of these subspaces, and is trivial
on the associated graded space. Similarly, the vector
space $\P\ts(\CC^{\ts m+l}\ot\CC^{\ts n})$ admits an ascending
filtration by $\q$-submodules such that $\q$ acts trivially on
each of the corresponding graded subspaces.
The latter filtration induces a filtration of $W$
by $\q$-submodules such that the corresponding graded quotient
$\operatorname{gr}W$ is a free $\U(\q)$-module. The
space of coinvariants $(\ts\operatorname{gr}W)_{\ts\q}$ is therefore
isomorphic to $V\ns\ot U\ot\P\ts(\CC^{\ts m+l}\ot\CC^{\ts n})\ts$,
via the bijective linear map
$$
y\ot x\ot f\,\mapsto\,1\ot y\ot x\ot f\,+\,\q\cdot(\operatorname{gr}W)\ts.
$$
Therefore the linear map $\chi$ is bijective as well.

Let us now demonstrate that the map $\chi$ intertwines the actions of
the Yangian $\Y(\gl_n)\ts$. Consider the $(m+l)\times(m+l)$ matrix whose
$ab\ts$-entry is $\de_{ab}\ts u+\ns E_{\ts ba}\ts$. Here we regard 
$E_{ba}$ as an element of the algebra $\U(\gl_{m+l})\ts$. 
Write this matrix in the block form 
\eqref{blockmat} where $A,B,C,D$ are matrices of sizes
$m\times m$, $m\times l$, $l\times m$, $l\times l$ respectively.
In the notation introduced in the end of Section 1, 
here $A=u+\Ep$.  
Using the observation made there along with the definition \eqref{1.33}
of the comultiplication, the action of the algebra $\Y(\gl_N)$ on
the vector space \eqref{vuprod} of the tensor product of the two
$\Y(\gl_n)$-modules $\A_{\ts m}(V)$ and $\A_{\ts l}^{\ts m}\ts(U)$
can be described by assigning to every series $T_{ij}(u)$ the 
product of the series 
$$
\sum_{k=1}^n\,\,
\Bigl(\ts\de_{ik}\,+\!
\sum_{a,b=1}^m
(A^{-1})_{ab}\ot x_{ai}\,\d_{\ts bk}\ts\Bigr)
\Bigl(\ts\de_{kj}\,+\!
\sum_{c,d=1}^{l}
((\ts D-m\ts )^{-1})_{cd}\ot x_{\ts m+c,k}\,\d_{\ts m+d,j}\ts\Bigr)
$$
\begin{equation}
\label{sum1}
=\,\de_{ij}\,+\!\ts
\sum_{a,b=1}^m
(A^{-1})_{ab}\ot x_{ai}\,\d_{\ts bj}\,+
\sum_{c,d=1}^{l}
((\ts D-m\ts )^{-1})_{cd}\ot x_{\ts m+c,i}\,\d_{\ts m+d,j}\ +
\vspace{4pt}
\end{equation}
\begin{equation}
\label{sum2}
\sum_{k=1}^n\ \,
\sum_{a,b=1}^m\,
\sum_{c,d=1}^{l}\ \,
(A^{-1})_{ab}\,((\ts D-m\ts )^{-1})_{cd}\ot 
x_{ai}\,\d_{\ts bk}\,x_{\ts m+c,k}\,\d_{\ts m+d,j}\,.
\end{equation}
Note that in \eqref{sum2} we have
$\d_{\ts bk}\,x_{\ts m+c,k}=x_{\ts m+c,k}\,\d_{\ts bk}$ because $b\le m\ts$.
The first tensor factors of all summands in
\eqref{sum1} and \eqref{sum2} correspond to the
action of the universal enveloping algebra
$\U(\ts\gl_m\op\gl_{\ts l}\ts)$ on $V\ns\ot U$.

Let us now write the matrix inverse to \eqref{blockmat} as
the block matrix
$$
\begin{bmatrix}
\ \At\,&\Bt\ 
\\
\ \Ct\,&\Dt\ 
\end{bmatrix}
$$
where $\At,\Bt,\Ct,\Dt$ are matrices of sizes
$m\times m$, $m\times l$, $l\times m$, $l\times l$ respectively.
Each of these four blocks is regarded 
as formal power series in $u^{-1}$ 
with matrix coefficients. The entries of these matrix coefficients 
belong to the algebra $\U(\ts\gl_{m+l})\ts$.
By once again using the observation made in the end
of Section 1, the action of $\Y(\gl_n)$
on the vector space $W$ can now be described by assigning to
every series $T_{ij}(u)$ the sum of the series 
$$
\de_{ij}+\sum_{a,b=1}^m
\At_{\ts ab}\ot x_{ai}\,\d_{\ts bj}\,+\,
\sum_{a=1}^m\,\sum_{d=1}^{l}\, 
\Bt_{\ts ad}\ot x_{ai}\,\d_{\ts m+d,j}\ +
$$
$$
\sum_{b=1}^m\,\sum_{c=1}^{l}\, 
\Ct_{\ts cb}\ot x_{\ts m+c,i}\,\d_{\ts bj}\,+
\sum_{c,d=1}^{l}
\Dt_{\ts cd}\ot x_{\ts m+c,i}\,\d_{\ts m+d,j}\,.
$$
The first tensor factors in all summands here correspond to the
action of the algebra $\U(\gl_{m+l})$ on the vector space
$\U(\q)\ot V\ns\ot U$ of the parabolically induced module $V\bt U\ts$.

Let us apply these tensor factors to elements of the
subspace \eqref{onesub}. By Lemma \ref{lemma1}, %we have
$$
\At=(A-B\ts D^{-1}\ts C)^{-1}\ts.
$$
All entries of the matrix $C$ belong to $\qp$ and
hence act on the subspace \eqref{onesub} as zeroes.
Further, we have $A=u+\Ep$. Every entry of the matrix $\Ep$
belongs to the subalgebra $\gl_m\subset\gl_{m+l}$ and
the adjoint action of this subalgebra on $\gl_{m+l}$
preserves $\qp$. Therefore the results of applying
$(A^{-1})_{ab}$ and $\At_{ab}$ to elements of
the subspace \eqref{onesub} are the same.
Similar arguments show that any entry of the matrix
$$
\Ct=-\ts D^{-1}C\ts(A-B\ts D^{-1}\ts C)^{-1}=-\ts D^{-1}C\ts\At
$$
act on the subspace \eqref{onesub} as zero.

Consider the matrix 
$$
\Dt=(\ts D-C\ts A^{-1}B)^{-1}\ts.
$$
In the notation of Lemma \ref{lemma2} the $ab\ts$-entry of the matrix
$A^{-1}$ is $X_{ab}(u)\ts$, and the trace of $A^{-1}$ is $Z(u)\ts$. 
Using that lemma, the $cd\ts$-entry of the $l\times l$ 
matrix $D-C\ts A^{-1}B$ equals
$$
\de_{cd}\,u\ts+E_{m+d,m+c}-\sum_{a,b=1}^m 
E_{\ts a,m+c}\,X_{ab}(u)\,E_{\ts m+d,b}\,=
\de_{cd}\,u\ts+E_{m+d,m+c}
\vspace{-4pt}
$$
$$
-\sum_{a,b=1}^m 
E_{\ts a,m+c}\,E_{\ts m+d,b}\,X_{ab}(u)\,(\ts1-Z(u))^{-1}=\,
\de_{cd}\,u\ts+E_{m+d,m+c}
$$
$$
-\sum_{a,b=1}^m 
(\ts E_{\ts m+d,b}\,E_{\ts a,m+c}+
\de_{cd}\,E_{ab}-\de_{ab}\,E_{\ts m+d,m+c}\ts)
\,X_{ab}(u)\,(\ts1-Z(u))^{-1}=
\vspace{7pt}
$$
\begin{equation}
\label{dmcd}
\de_{cd}\,(\ts u-(\ts m-u\,Z(u))\ts(\ts1-Z(u))^{-1})\,+
E_{m+d,m+c}\,(\ts1+Z(u)\ts(\ts1-Z(u))^{-1})
\vspace{12pt}
\end{equation}
\begin{equation}
\label{efactor}
-
\sum_{a,b=1}^m 
E_{\ts m+d,b}\,E_{\ts a,m+c}
\,X_{ab}(u)\,(\ts1-Z(u))^{-1}.
\vspace{10pt}
\end{equation}
We used the identity
$$
\sum_{a,b=1}^m E_{ab}\,X_{ab}(u)\,=\,m-u\,Z(u)
$$
which follows from the definitions \eqref{xabu} and \eqref{zu}.
The expression in the line \eqref{dmcd} is equal to
$$
(\ts D-m)_{\ts cd}\ts(\ts1-Z(u))^{-1}.
$$
The factor $E_{\ts a,m+c}$ in any summand in the line \eqref{efactor}
belongs to $\qp$ while every element of $\qp$
acts on the subspace \eqref{onesub} as zero. 
The coefficients of the series $X_{ab}(u)$ and $Z(u)$ in \eqref{efactor}
belong to $\U(\gl_m)$ while the adjoint action of subalgebra 
$\gl_m\subset\gl_{m+l}$ preserves $\qp$. The adjoint action of
every element $E_{m+d,m+c}\in\gl_{m+l}$ also preserves $\qp$.
Therefore the result of applying $\Dt_{cd}$ to elements of
the subspace \eqref{onesub} is the same as that of applying
$$
(\ts1-Z(u))\ts((\ts D-m)^{-1})_{cd}\ts.
$$

Now consider
$$
\Bt=-\ts A^{-1}B\ts(\ts D-C\ts A^{-1}B)^{-1}=-\ts A^{-1}B\ts\Dt\ts.
$$
The above arguments show that the result of applying
the $ad\ts$-entry of this matrix to
elements of the subspace \eqref{onesub} is the same as that of applying
the $ad\ts$-entry of the matrix
$$
-\ts A^{-1}B\ts(\ts1-Z(u))\ts(\ts D-m)^{-1}\ts.
$$
By using Lemma \ref{lemma2} once again, the latter entry equals
$$
-\,\sum_{b=1}^m\,\sum_{c=1}^l\,
X_{ab}(u)\,E_{\ts m+c,b}
(\ts1-Z(u))\ts((\ts D-m)^{-1})_{cd}\,=
$$
$$
-\,\sum_{b=1}^m\,\sum_{c=1}^l\,
E_{\ts m+c,b}\,X_{ab}(u)\ts((\ts D-m)^{-1})_{cd}\,.
$$

Thus we have proved that
the action of $\Y(\gl_n)$ on the elements of the subspace
\begin{equation}
\label{1vu}
1\ot V\ns\ot U\ot\P\ts(\CC^{\ts m+l}\ot\CC^{\ts n})\subset W
\end{equation}
can be described by assigning to every series $T_{ij}(u)$ 
the sum of the series 
$$
\de_{ij}+\sum_{a,b=1}^m
(A^{-1})_{\ts ab}\ot x_{ai}\,\d_{\ts bj}\,+
\sum_{c,d=1}^{l}
(1-Z(u))\ts((\ts D-m)^{-1})_{\ts cd}\ot x_{\ts m+c,i}\,\d_{\ts m+d,j}
\vspace{-6pt}
$$
\begin{equation}
\label{exud}
-\,\sum_{a,b=1}^m\,\sum_{c,d=1}^{l}\, 
E_{\ts m+c,b}\,X_{ab}(u)\ts((\ts D-m)^{-1})_{cd}
\ot x_{ai}\,\d_{\ts m+d,j}\ts.
\vspace{4pt}
\end{equation}
Let us now consider the results of the action of $\Y(\gl_n)$
on this subspace modulo $\q\cdot W$. Since $E_{\ts m+c,b}\in\q\,$,
the expession
in the line \eqref{exud} can be then replaced by
$$
\sum_{k=1}^n\ \sum_{a,b=1}^m\,\sum_{c,d=1}^{l}\, 
X_{ab}(u)\ts((\ts D-m)^{-1})_{cd}
\ot x_{m+c,k}\,\d_{\ts bk}\,x_{ai}\,\d_{\ts m+d,j}\,=
$$
$$
\sum_{k=1}^n\ \sum_{a,b=1}^m\,\sum_{c,d=1}^{l}\, 
(A^{-1})_{ab}\ts((\ts D-m)^{-1})_{cd}
\ot x_{m+c,k}\,x_{ai}\,\d_{\ts bk}\,\d_{\ts m+d,j}\ +
$$
$$
\sum_{c,d=1}^{l}\, 
Z(u)\ts((\ts D-m)^{-1})_{cd}
\ot x_{m+c,i}\,\d_{\ts m+d,j}\ts.
$$
Here we used the equality of differential operators
$\d_{\ts bk}\,x_{ai}=
x_{ai}\,\d_{\ts bk}+\de_{ab}\,\de_{ik}\ts$.
By making this replacement we show 
that modulo $\q\cdot W$,
the action of $\Y(\gl_n)$ on elements of
the subspace \eqref{1vu}
can be described by assigning to every series $T_{ij}(u)$ 
the sum of the series 
$$
\de_{ij}+\sum_{a,b=1}^m
(A^{-1})_{\ts ab}\ot x_{ai}\,\d_{\ts bj}\,+
\sum_{c,d=1}^{l}
(1-Z(u))\ts((\ts D-m)^{-1})_{\ts cd}\ot x_{\ts m+c,i}\,\d_{\ts m+d,j}\ +
\vspace{-2pt}
$$
$$
\sum_{k=1}^n\ \sum_{a,b=1}^m\,\sum_{c,d=1}^{l}\, 
(A^{-1})_{ab}\ts((\ts D-m)^{-1})_{cd}
\ot x_{m+c,k}\,x_{ai}\,\d_{\ts bk}\,\d_{\ts m+d,j}\ +
$$
$$
\sum_{c,d=1}^{l}\, 
Z(u)\ts((\ts D-m)^{-1})_{cd}
\ot x_{m+c,i}\,\d_{\ts m+d,j}
$$
which is equal to the sum of the series in the lines 
\eqref{sum1} and \eqref{sum2}.
This equality proves that the map $\chi$ intertwines the actions of
the Yangian $\Y(\gl_n)\ts$.
\qed

\medskip\smallskip
By the transitivity of induction, Theorem \ref{parind}
can be extended from the maximal to all
parabolic subalgebras of the Lie algebra $\gl_m\ts$.
Consider the case of the Borel
subalgebra $\h\op\np$ of $\gl_m\ts$.
Here $\h$ is the Cartan subalgebra of $\gl_m$ spanned
by the elements $E_{aa}\ts$, whereas
$\n'$ is the nilpotent subalgebra of $\gl_m$ spanned
by the elements $E_{ab}$ with $a<b\ts$.

Take any element $\mu$ of the vector space $\h^\ast$ dual to
$\h\ts$, any such element is called a \textit{weight}\ts.
The weight $\mu$ can be identified with
the sequence $(\ts\mu_1\lcd\mu_m)$ of its \textit{labels},
where $\mu_a=\mu(E_{aa})$ for each $a=1\lcd m\ts$.
Consider the Verma module $M_{\ts\mu}$ over the Lie algebra
$\gl_m\ts$. It can be described as the quotient of the algebra
$\U(\gl_m)$ by the left ideal generated by all the elements
$E_{ab}$ with $a<b$ and the elements $E_{aa}-\mu_a\ts$.
The elements of the Lie algebra $\gl_m$ act on this quotient
via left multiplication. The image of the element $1\in\U(\gl_m)$
in this quotient is denoted by $1_\mu\ts$. Then
$X\cdot 1_\mu=0$ for all $X\in\n'$ while
$$
X\cdot 1_\mu=\mu\ts(X)\cdot 1_\mu
\quad\textrm{for all}\quad
X\in\h\ts.
$$

Let us 
apply the functor \eqref{zelefun} to the $\gl_m$-module $V=M_{\ts\mu}\ts$,
and the functor \eqref{drinfun} to the resulting $\H_N$-module
$$
W=(\ts M_{\ts\mu}\ot(\CC^{\ts m})^{\ot N})_{\ts\n}\ts.
$$
We obtain the $\Y(\gl_n)$-module
$$
((\ts M_{\ts\mu}\ot(\CC^{\ts m})^{\ot N})_{\ts\n}\ot
(\CC^{\ts n})^{\ot N})^{\Sym_N}
=
(\ts M_{\ts\mu}\ot\S^N(\CC^{\ts m}\ot\CC^{\ts n}))_{\ts\n}\ts.
$$
By taking the direct sum over $N=0,1,2,\ts\ldots$
of these $\Y(\gl_n)$-modules, we obtain the $\Y(\gl_n)$-module
$$
\A_{\ts m}(\ts M_{\ts\mu})_{\ts\n}=
(\ts M_{\ts\mu}\ot\S\ts(\CC^{\ts m}\ot\CC^{\ts n}))_{\ts\n}\,.
$$
Note that $\A_{\ts m}(\ts M_{\ts\mu})_{\ts\n}$ is also a module over
the Cartan subalgebra $\h\subset\gl_m\ts$. Using the basis
$E_{11}\lcd E_{mm}$ identify $\h$ with the direct sum of $m$
copies of the Lie algebra $\gl_{\ts1}\ts$.
Consider the Verma modules $M_{\ts\mu_1}\lcd M_{\ts\mu_m}$ over
$\gl_{\ts1}\ts$. By applying
Theorem~\ref{parind} repeatedly we get the next result,
which can be also derived from \cite[Theorem 3.3.1]{AS}.

\begin{corollary}
\label{bimequiv}
The bimodule $\A_{\ts m}(\ts M_{\ts\mu})_{\ts\n}$ of\/ $\h$ and\/ $\Y(\gl_n)$
is equivalent to the tensor product
$$
\A_{\ts1}\ts(\ts M_{\ts\mu_1})
\ot
\A_{\,1}^{\tts1}\ts(\ts M_{\ts\mu_2})
\ot
\ldots
\ot
\A_{\,1}^{\,m-1}\ts(\ts M_{\ts\mu_m})\ts.
$$ 
\end{corollary}

We complete this section with describing 
for any $t\com z\in\CC$ the bimodule $\A_{\ts1}^{\ts z}\ts(M_{\ts t})$ 
over $\gl_{\ts1}$ and $\Y(\gl_n)\ts$.
The Verma module $M_{\ts t}$ over
$\gl_{\ts1}$ is one-dimensional, and the element
$E_{11}\in\gl_{\ts1}$ acts on $M_{\ts t}$ as multiplication by $t\ts$.
The vector space of bimodule 
$\A_{\ts1}(M_{\ts t})$ is the symmetric algebra 
$\S\ts(\CC^{\ts1}\ot\CC^{\ts n})=\S\ts(\CC^{\ts n})\ts$,
which we identify with 
$\P\ts(\CC^{\ts1}\ot\CC^{\ts n})=\P\ts(\CC^{\ts n})\ts$. 
Then $E_{11}$ acts on $\A_{\ts1}(M_{\ts t})$ as the differential operator
$$
t\,+\,\sum_{k=1}^n\,x_{1k}\,\d_{\ts1k}\,.
$$
The action of $E_{11}$ on 
$\A_{\ts1}^{\ts z}\ts(M_{\ts t})$ is the same as on $\A_{\ts1}(M_{\ts t})\ts$.
The generator $T_{ij}^{\ts(s+1)}$ of $\Y(\gl_n)$
with $s=0,1,2,\ts\ldots$
acts on $\A_{\ts1}(M_{\ts t})$ as the differential operator 
$$
(-t)^{\ts s}\,x_{1i}\,\d_{\ts1j}\,,
$$
this is what Proposition \ref{dast} states in the case $m=1$.
Note that the operator $x_{1i}\,\d_{\ts1j}$ describes the
action on $\P\ts(\CC^{\ts1}\ot\CC^{\ts n})=\P\ts(\CC^{\ts n})$ 
of the element $E_{ij}\in\gl_n\ts$.
Hence the action of the algebra $\Y(\gl_n)$ on 
$\A_{\ts1}(M_{\ts t})$ can be obtained from the action
of $\gl_n$ on $\P\ts(\CC^{\ts n})\ts$ by pulling back through 
the homomorphism $\pi_n:\Y(\gl_n)\to\U(\gl_n)\ts$, and then
through the automorphism $\tau_{-t}$ of $\Y(\gl_n)\ts$; see
the definitions \eqref{tauz} and \eqref{pin}.
Hence the action of $\Y(\gl_n)$ on $\A_{\ts1}^{\ts z}\ts(M_{\ts t})$
can be obtained from the action
of $\gl_n$ on $\P\ts(\CC^{\ts n})\ts$ by pulling back through 
$\pi_n$ and then through the automorphism $\tau_{z-t}\ts$.

%------------------------------------------------------------------------------

\section*{\normalsize 3.\ Zhelobenko operators}
\setcounter{section}{3}
\setcounter{equation}{0}
\setcounter{theorem}{0}

Consider the symmetric group $\Sym_m$ as the Weyl group of
the reductive Lie algebra $\gl_m\ts$. This group acts on the 
vector space $\gl_m$ so that for any $\si\in\Sym_m$ and $a,b=1\lcd m$
$$
\si:E_{ab}\mapsto E_{\si(a)\ts\si(b)}\ts.
$$
This action extends to an action of the group $\Sym_m$
by automorphisms of the associative algebra $\U(\gl_m)\ts$.
The group $\Sym_m$ also acts on the vector space $\h^\ast\ts$.
Let $E_{11}^{\,\ast}\lcd E_{mm}^{\,\ast}$ be the basis
of $\h^\ast$ dual to the basis $E_{11}\lcd E_{mm}$ of $\h\ts$. Then
$$
\si:E_{aa}^{\,\ast}\mapsto E_{\si(a)\ts\si(a)}^{\,\ast}\ts.
$$
If we identify each weight $\mu\in\h^\ast$ with
the sequence $(\ts\mu_1\lcd\mu_m)$ of its labels, then
$$
\si:(\ts\mu_1\lcd\mu_m)\mapsto
(\,\mu_{\ts\si^{-1}(1)}\lcd\mu_{\ts\si^{-1}(m)}\ts)\ts.
$$
Let $\rho\in\h^\ast$ be the weight with
sequence of labels $(\ts0,-1\lcd1-m)\ts$.
The \textit{shifted\/} action of any element $\si\in\Sym_m$ on $\h^\ast$
is defined by the assignment
\begin{equation}
\label{saction}
\mu\,\mapsto\,\si\circ\mu=\si\ts(\mu+\rho)-\rho\ts.
\end{equation}

For any $a,b=1\lcd m$ put $\ep_{ab}=E_{aa}^{\,\ast}-E_{bb}^{\,\ast}\,$.
The elements $\ep_{ab}\in\h^\ast$ with $a<b$ and $a>b$ are
the \textit{positive\/} and \textit{negative roots\/} 
respectively. Note that $\ep_{ab}=0$ when $a=b\ts$.
The elements $\ep_c=\ep_{c,c+1}\in\h^\ast$ with $c=1\lcd m-1$ are the
\textit{simple\/} positive roots. Put
$$
E_{c}=E_{c,c+1}\ts,
\quad
F_{c}=E_{c+1,c}
\quad\text{and}\quad
H_c=E_{cc}-E_{c+1,c+1}\ts.
$$
For any $a=1\lcd m-1$ these three elements of the
Lie algebra $\gl_m$ span a subalgebra 
isomorphic to the Lie algebra $\mathfrak{sl}_{\ts2}\ts$.

For any $\gl_m$-module $V$ and any $\la\in\h^\ast$
a vector $v\in V$ is said to be \textit{of weight} $\la$ if
$X\,v=\la\ts(X)\,v$ for any $X\in\h\ts$. We will
denote by $V^\la$ the subspace in $V$ formed by all
vectors of weight $\la\ts$. Recall that $\n$ denotes
the nilpotent subalgebra of $\gl_m$ spanned
by the elements $E_{ab}$ with $a>b\ts$.  
In this section, we will employ
the general notion of a \textrm{Mickelsson algebra} introduced in \cite{M1}
and developed by D.\ts Zhelobenko \cite{Z}. 
Namely, we will show how this notion gives rise to
a distinguished $\Y(\gl_n)$-intertwining operator
\begin{equation}
\label{distoper}
\A_{\ts m}(\ts M_{\ts\mu})_{\ts\n}^{\ts\la}
\ \to\,
\A_{\ts m}(\ts M_{\ts\si\,\circ\ts\mu})_{\ts\n}^{\,\si\,\circ\ts\la}
\end{equation}
for any element $\si\in\Sym_m$ and any
weight $\mu\in\h^\ast$ such that 
\begin{equation}
\label{muass}
\mu_a-\mu_b\notin\ZZ
\quad\text{whenever}\quad
a\neq b\ts.
\end{equation}
Note that the source and the target vector
spaces in \eqref{distoper} are non-zero only if
all labels of the weight $\la-\mu$ are non-negative integers. Then 
$\la_a-\la_b\notin\ZZ$ whenever $a\neq b\ts$.

We have a representation
$\ga:\,\U(\gl_m)\to\PD\ts(\CC^{\ts m}\ot\CC^{\ts n})$
such that the image $\ga(E_{ab})$ is the differential operator \eqref{glmact}.
Note that the group $\Sym_m$ acts by automorphisms of the algebra 
$\PD\ts(\CC^{\ts m}\ot\CC^{\ts n})\ts$, so that for $k=1\lcd n$
$$
\si:\
x_{ak}\ts\mapsto\ts x_{\si(a)\ts k}\,,\
\d_{\ts bk}\ts\mapsto\d_{\ts\si(b)\ts k}\,.
$$
The homomorphism $\ga$ is $\Sym_m$-equivariant. Let $\Ar$ be
the associative algebra generated by the algebras
$\U(\gl_m)$ and $\PD\ts(\CC^{\ts m}\ot\CC^{\ts n})$
with the cross relations
\begin{equation}
\label{defar}
[X\com Y]=[\ts\ga(X)\com Y]
\end{equation}
for any $X\in\gl_m$ and $Y\in\PD\ts(\CC^{\ts m}\ot\CC^{\ts n})\ts$.
Here each pair of the square brackets denotes the commutator in $\Ar\ts$.
Note that the algebra $\Ar$ is isomorphic to the 
the tensor product of associative algebras \eqref{tenprod}.
The isomorphism can be defined by mapping
the elements $X\in\gl_m$ and $Y\in\PD\ts(\CC^{\ts m}\ot\CC^{\ts n})$
of the algebra $\Ar$ respectively to the elements 
$$
X\ot1+1\ot\ga(X)
\quad\textrm{and}\quad
1\ot Y
$$
of the algebra \eqref{tenprod}.
This isomorphism is $\Sym_m$-equivariant, and the image 
of the element $E_{ab}\in\gl_m$ under this isomorphism
equals \eqref{eabact}. We will use this isomorphism later on.

Let $\J$ be the right ideal of the algebra $\Ar$
generated by the elements of the subalgebra $\n\subset\gl_m\ts$.
Let ${\rm Norm}\ts(\J)\ts\subset\ts\Ar$
be the normalizer of this right ideal, so that
$Y\in{\rm Norm}\ts(\J)$ if and only if $\,Y\ts\J\subset\J\ts$.
Then $\J$ is a two-sided ideal of ${\rm Norm}\ts(\J)\ts$.
Our particular \textit{Mickelsson algebra\/} is the quotient 
\begin{equation}
\label{malg}
\R\,=\,\J\,\backslash\,{\rm Norm}\ts(\J)\ts.
\end{equation}

% Note that ${\rm Norm}\ts(\J)$ contains any element of the centralizer
% of the subalgebra $\U(\gl_m)$ of $\Ar\ts$.
% In particular, for any $s=0,1,2,\ts\ldots$ the normalizer 
% ${\rm Norm}\ts(\J)$ contains the
% element \eqref{combact}, which describes the action
% of the generator $T_{ij}^{(s+1)}$ of the Yangian $\Y(\gl_n)$ 
% on the vector space 
% $V\ot\ts\P\ts(\CC^{\ts m}\ot\CC^{\ts n})\ts$.

\medskip\smallskip
\noindent\textit{Remark.}
Via its isomorphism with \eqref{tenprod},
the associative algebra $\Ar$ acts on the tensor product
$V\ot\ts\P\ts(\CC^{\ts m}\ot\CC^{\ts n})$
for any $\gl_m$-module $V$.
The defining embedding of $\gl_m$ into $\Ar$ corresponds to the diagonal action
of the Lie algebra $\gl_m$ on this tensor product.
The algebra $\R$ then acts
on the space of $\n$-coinvariants of $\gl_m$-module 
$V\ot\ts\P\ts(\CC^{\ts m}\ot\CC^{\ts n})\ts$.
\qed

\medskip\smallskip
Let $\,\overline{\!\U(\h)\!\!\!}\,\,\,$ 
be the ring of fractions of the %commutative 
algebra $\U(\h)$ relative to the denominators set
\begin{equation}
\label{denset}
\{\,E_{aa}-E_{bb}+z\ |\ 1\le a\com b\le m\ts;\ a\neq b\ts;\ z\in\ZZ\,\ts\}\,.
\end{equation}
The elements of this ring can also
be regarded as rational functions on the vector space
$\h^\ast\ts$. The elements of 
$\U(\h)\subset\,\overline{\!\U(\h)\!\!\!}\,\,\,$
are then regarded as polynomial functions on $\h^\ast\ts$.
Denote by $\Ab$ the ring of fractions of $\Ar$ 
relative to the same set of denominators \eqref{denset},
regarded as elements of $\Ar$ using the embedding of $\h\subset\gl_m$ 
into $\Ar\ts$. The ring $\Ab$ is defined 
due to the following relations
in $\U(\gl_m)$ and $\Ar\ts$: for $a,b=1\lcd m$
and any $H\in\h$
$$
[\ts H\ts,E_{ab}\ts]=\ep_{ab}(H)\ts E_{ab}\,,
\quad\
[\ts H\ts,x_{ak}\ts]=
E_{aa}^{\,\ast}(H)\,x_{ak}\,,
\quad
[\ts H\ts,\d_{\ts bk}\ts]=
-\ts E_{\ts bb}^{\,\ast}(H)\,\d_{\ts bk}\,.
$$
Therefore the ring $\Ar$ satisfies the Ore condition relative
to its subset \eqref{denset}.
Using left multiplication by elements of
$\,\overline{\!\U(\h)\!\!\!}\,\,\,$,
the ring of fractions $\Ab$ becomes a module over
$\,\overline{\!\U(\h)\!\!\!}\,\,\,\ts$.

The ring $\Ab$ is also an associative algebra over the field $\CC\ts$.
For each $c=1\lcd m-1$ define a linear map $\xi_{\ts c}:\Ar\to\Ab$
by setting
\begin{equation}
\label{q1}
\xi_{\ts c}\,(\ts Y)=Y+\,\sum_{s=1}^\infty\,\,
(\ts s\ts!\,H_c^{\ts(s)}\ts)^{-1}\ts E_c^{\ts s}\,
\widehat{F}_c^{\ts s}(\ts Y)
\end{equation}
for $Y\in\Ar\ts$. Here
$$
H_c^{\ts(s)}=H_c(H_c-1)\ldots(H_c-s+1)
$$
and $\widehat{F}_c$ is the operator of adjoint action
corresponding to the element $F_c\in\Ar\ts$, so that
$$
\widehat{F}_c(\ts Y)=[\ts F_c\ts\com Y\ts]\ts.
$$
For any given element $Y\in\Ar$ only finitely many terms of the sum 
\eqref{q1} differ from zero, hence the map $\xi_{\ts c}$ is well defined.
The definition \eqref{q1} and the following proposition
go back to \cite[Section 2\ts]{Z}.
Denote $\Jb=\,\overline{\!\U(\h)\!\!\!}\,\,\,\,\J\ts$.
Then $\Jb$ is a right ideal of the algebra $\Ab\ts$.

\begin{proposition}
\label{prop1}
For any $X\in\h$ and $Y\in\Ar$ we have
\begin{equation}
\label{q11}
\xi_{\ts c}(X\ts Y)\in
(\ts X+\ep_c(X))\,\ts\xi_{\ts c}(\ts Y)\ts+\ts\Jb\ts,
\end{equation}
\begin{equation}
\label{q12}
\xi_{\ts c}(\ts Y X)\in\,
\xi_{\ts c}(\ts Y)\ts(\ts X+\ep_c(X))\ts+\ts\Jb\ts.
\end{equation}
\end{proposition}

\begin{proof}
It suffices to verify each of the properties \eqref{q11} and \eqref{q12}
for two elements of $\h\ts$, 
$$
X=E_{cc}+E_{c+1,c+1}
\ \quad\text{and}\ \quad
X=E_{cc}-E_{c+1,c+1}=H_c\ts.
$$
For the first of these two elements we have
$[\ts E_c\,,X\ts]=[\ts F_c\,,X\ts]=0$ and $\ep_c(X)=0\ts$,
so that the properties \eqref{q11} and \eqref{q12} are obvious. 

For $X=H_c$ the proof of \eqref{q11}
is based on the following commutation relations in the subalgebra of 
$\U(\gl_m)$ generated by the three elements $E_c,F_c$ and $H_c\ts$:
for $s=1,2,\ts\ldots$
\begin{eqnarray}
\label{had}
[\,E_c^{\ts s}\ts,\ns H_c\ts]
&=&
-2\ts s\,E_c^{\ts s}\,,
\\
\label{fad}
[\,E_c^{\ts s}\ts,F_c\,]
&=&
s\,(H_c-s+1)\,E_c^{\ts s-1}\ts.
\end{eqnarray}

Let us use the symbol $\ts\,\equiv\,$ to indicate equalities in the algebra
$\Ab$ modulo the ideal $\Jb\ts$. By the definition \eqref{q1}, for
any element $Y\in\Ar$ we have
$$
\begin{aligned}
\xi_{\ts c}\,(\ts H_c\,Y)
&\,=\,H_c\,Y+\,
\sum_{s=1}^\infty\,\,
(\ts s\ts!\,H_c^{\ts(s)}\ts)^{-1}\ts E_c^{\ts s}\,
\widehat{F}_c^{\ts s}(\ts H_c\,Y)
\\
&\,=\,H_c\,Y+\,
\sum_{s=1}^\infty\,\,
(\ts s\ts!\,H_c^{\ts(s)}\ts)^{-1}\ts E_c^{\ts s}\,
(\ts H_c\,\widehat{F}_c^{\ts s}(\ts Y)+
2\ts s\,F_c\,\widehat{F}_c^{\ts s-1}(\ts Y))
\\
&\,\equiv\,H_c\,Y+\,
\sum_{s=1}^\infty\,\,
(\ts s\ts!\,H_c^{\ts(s)}\ts)^{-1}\ts
((\ts H_c-2\ts s)\,E_c^{\ts s}\,\widehat{F}_c^{\ts s}(\ts Y)
\\
&\phantom{\,=\,H_c\,Y}\ +\,
\sum_{s=1}^\infty\,\,
(\ts s\ts!\,H_c^{\ts(s)}\ts)^{-1}\ts
\cdot
2\ts s^2\,(H_c-s+1)\,E_c^{\ts s-1}\,\widehat{F}_c^{\ts s-1}(\ts Y))
\\
&\,=\,H_c\,Y+\,
\sum_{s=1}^\infty\,\,
(\ts s\ts!\,H_c^{\ts(s)}\ts)^{-1}\ts
((\ts H_c-2\ts s)\,E_c^{\ts s}\,\widehat{F}_c^{\ts s}(\ts Y)
\\
&\phantom{\,=\,H_c\,Y}\ +\,
\sum_{s=0}^\infty\,\,
(\ts s\ts!\,H_c^{\ts(s)}\ts)^{-1}\ts
\cdot
2\,(s+1)\,E_c^{\ts s}\,\widehat{F}_c^{\ts s}(\ts Y))
\\
&\,=\,(\ts H_c+2\,)\,\ts\xi_{\ts c}\,(\ts Y)
\,=\,(\ts H_c+\ep_c(H_c))\,\ts\xi_{\ts c}\,(\ts Y)
\phantom{\sum^\infty}
\end{aligned}
$$
as needed. To get the equivalence relation above, we also used the inclusion
$\ts\,\overline{\!\U(\h)\!\!\!}\,\,\,\ts F_c\ts\subset\ts\Jb\ts$.

Similarly, for any element $Y\in\Ar$ we have
$$
\begin{aligned}
\xi_{\ts c}\,(\,YH_c\ts)
&\,=\,YH_c\ts+\,
\sum_{s=1}^\infty\,\,
(\ts s\ts!\,H_c^{\ts(s)}\ts)^{-1}\ts E_c^{\ts s}\,
\widehat{F}_c^{\ts s}(\,YH_c\ts)
\\
&\,=\,YH_c\ts+\,
\sum_{s=1}^\infty\,\,
(\ts s\ts!\,H_c^{\ts(s)}\ts)^{-1}\ts E_c^{\ts s}\,
(\,\widehat{F}_c^{\ts s}(\ts Y)\ts H_c\ts+
2\ts s\,\widehat{F}_c^{\ts s-1}(\ts Y)\ts F_c\ts)
\\
&\,=\,YH_c\ts+\,
\sum_{s=1}^\infty\,\,
(\ts s\ts!\,H_c^{\ts(s)}\ts)^{-1}\ts E_c^{\ts s}\,
(\,\widehat{F}_c^{\ts s}(\ts Y)\ts H_c\ts-
2\ts s\,\widehat{F}_c^{\ts s}(\ts Y)+
2\ts s\,F_c\,\widehat{F}_c^{\ts s-1}(\ts Y))
\\
&\,\equiv\,YH_c\ts+\,
\sum_{s=1}^\infty\,\,
(\ts s\ts!\,H_c^{\ts(s)}\ts)^{-1}\ts E_c^{\ts s}\,
(\,\widehat{F}_c^{\ts s}(\ts Y)\ts H_c\ts-
2\ts s\,\widehat{F}_c^{\ts s}(\ts Y))
\\
&\phantom{\,=\,YH_c}\ +\,
\sum_{s=1}^\infty\,\,
(\ts s\ts!\,H_c^{\ts(s)}\ts)^{-1}\cdot
2\ts  s^2\,(H_c-s+1)\,E_c^{\ts s-1}\,
\widehat{F}_c^{\ts s-1}(\ts Y))
\\
&\,=\,YH_c\ts+\,
\sum_{s=1}^\infty\,\,
(\ts s\ts!\,H_c^{\ts(s)}\ts)^{-1}\ts E_c^{\ts s}\,
(\,\widehat{F}_c^{\ts s}(\ts Y)\ts H_c\ts-
2\ts s\,\widehat{F}_c^{\ts s}(\ts Y))
\\
&\phantom{\,=\,YH_c}\ +\,
\sum_{s=0}^\infty\,\,
(\ts s\ts!\,H_c^{\ts(s)}\ts)^{-1}\cdot
2\,(s+1)\,E_c^{\ts s}\,
\widehat{F}_c^{\ts s}(\ts Y))
\\
&
\,=\,\xi_{\ts c}\,(\ts Y)\ts(\ts H_c+2\,)
\,=\,\xi_{\ts c}\,(\ts Y)\ts(\ts H_c+\ep_c(H_c))\,.
\phantom{\sum^\infty}
\end{aligned}
$$
\vspace{-32pt}
\par\ \qed
\end{proof}

\smallskip\medskip
The property \eqref{q11} allows us to define a linear map 
$
\bar\xi_{\ts c}:\Ab\to\Jb\,\ts\backslash\,\Ab
$
by setting
$$
\bar\xi_{\ts c}(X\,Y)=Z\,\xi_{\ts c}(\ts Y)\ts+\ts\Jb
\quad\ \text{for any}\quad
X\in\,\overline{\!\U(\h)\!\!\!}\,\,\,
\quad\text{and}\quad
Y\in\Ar\ts,
$$
where the element $Z\in\,\overline{\!\U(\h)\!\!\!}\,\,\,$
is defined by the equality
$$
Z(\mu)=X(\ts\mu+\ep_c)
\quad\ \text{for any}\quad
\mu\in\h^\ast
$$
when $X$ and $Z$ are
regarded as rational functions on $\h^\ast\ts$.

The action of the group $\Sym_m$ on the algebra $\U(\h)$ 
extends to an action on $\,\overline{\!\U(\h)\!\!\!}\,\,\,\ts$, so that
for any $\si\in\Sym_m$
$$
(\ts\si\ts X)(\mu)=X(\ts\si^{-1}(\mu))
$$
when the element $X\in\,\overline{\!\U(\h)\!\!\!}\,\,\,$
is regarded as a rational function on $\h^\ast\ts$. The action
of $\Sym_m$ by automorphisms of the algebra $\Ar$ then 
extends to an action by automorphisms of $\Ab\ts$.
For any $c=1\lcd m-1$ let $\si_c\in\Sym_m$ be
the transposition of $c$ and $c+1\ts$.
Consider the image $\si_c(\ts\Jb\ts)\ts$,
this is again a right ideal of $\Ab$. Next proposition
also goes back to \cite{Z}. 

\begin{proposition}
\label{prop2}
We have\/ $\si_c(\ts\Jb\ts)\subset\ker\ts\bar\xi_{\ts c}\ts$.
\end{proposition}

\begin{proof}
Let $\n_{\ts c}$ be the vector subspace in $\gl_m$ spanned
by all the elements $E_{ab}$
where $a>b$ but $(a\com b)\neq(c+1\com c)\ts$. Then
$\si_c(\ts\Jb\ts)\subset\Ab$ is the right ideal generated
by the subspace $\n_{\ts c}\subset\gl_m\ts$ and by 
the element $E_c=E_{c,c+1}\ts.$ Here we use the embedding of
$\gl_m$ into $\Ab\ts$. Observe that
the subspace $\n_{\ts c}\subset\gl_m\ts$ is preserved by the
adjoint action of the elements $E_c,F_c$ and $H_c$ on $\gl_m\ts$.
Therefore $\xi_{\ts c}\ts(XY)\in\Jb$ for 
any $X\in\n_c$ and $Y\in\Ar\ts$, see \eqref{q1}. 
To prove Proposition \ref{prop2} it remains to show that 
$\xi_{\ts c}\ts(\ts E_c\ts Y\ts )\in\Jb$ for any $Y\in\Ar\ts$. We have
$$
\begin{aligned}
\xi_{\ts c}\,(\ts E_c\ts Y\ts)
&\,=\,E_c\ts Y\ts+\,
\sum_{s=1}^\infty\,\,
(\ts s\ts!\,H_c^{\ts(s)}\ts)^{-1}\ts E_c^{\ts s}\,
\widehat{F}_c^{\ts s}(\ts E_c\ts Y\ts)
\\
&\,=\,E_c\ts Y\ts+\,
\sum_{s=1}^\infty\,\,
(\ts s\ts!\,H_c^{\ts(s)}\ts)^{-1}\ts E_c^{\ts s}\ts
(\ts E_c\,\widehat{F}_c^{\ts s}(\ts Y\ts)-
s\,H_c\,\widehat{F}_c^{\ts s-1}(\ts Y\ts))
\\
&\phantom{\,=\,E_c\ts Y}\ -\,
\sum_{s=2}^\infty\,\,
(\ts s\ts!\,H_c^{\ts(s)}\ts)^{-1}\,E_c^{\ts s}\cdot
s\ts(s-1)\,F_c\,\widehat{F}_c^{\ts s-2}(\ts Y\ts))
\\
&\,\equiv\,E_c\ts Y\ts+\,
\sum_{s=1}^\infty\,\,
(\ts s\ts!\,H_c^{\ts(s)}\ts)^{-1}\ts 
(\,
E_c^{\ts s+1}\,\widehat{F}_c^{\ts s}(\ts Y\ts)-
s\ts(H_c-2\ts s)\,E_c^{\ts s}\,\widehat{F}_c^{\ts s-1}(\ts Y\ts))
\\
&\phantom{\,=\,E_c\ts Y}\ -\,
\sum_{s=2}^\infty\,\,
(\ts s\ts!\,H_c^{\ts(s)}\ts)^{-1}\ts 
s^2\ts(s-1)\ts(H_c-s+1)\,E_c^{\ts s-1}\,
\widehat{F}_c^{\ts s-2}(\ts Y\ts))
\\
&\,=\,E_c\ts Y\ts+\,
\sum_{s=1}^\infty\,\,
(\ts s\ts!\,H_c^{\ts(s)}\ts)^{-1}\ts 
(\,
E_c^{\ts s+1}\,\widehat{F}_c^{\ts s}(\ts Y\ts)-
s\ts(H_c-2\ts s)\,E_c^{\ts s}\,\widehat{F}_c^{\ts s-1}(\ts Y\ts))
\\
&\phantom{\,=\,E_c\ts Y}\ -\,
\sum_{s=1}^\infty\,\,
(\ts s\ts!\,H_c^{\ts(s)}\ts)^{-1}\ts 
s\ts(s+1)\,E_c^{\ts s}\,
\widehat{F}_c^{\ts s-1}(\ts Y\ts))
\\
&\,=\,E_c\ts Y\ts+\,
\sum_{s=1}^\infty\,\,
(\ts s\ts!\,H_c^{\ts(s)}\ts)^{-1}\ts 
(\,
E_c^{\ts s+1}\,\widehat{F}_c^{\ts s}(\ts Y\ts)-
s\ts(H_c-s+1)\,E_c^{\ts s}\,\widehat{F}_c^{\ts s-1}(\ts Y\ts))\,.
\end{aligned}
$$
The sum of all terms in the last line is equal to zero. 
Like in the proof of Proposition~\ref{prop1}, here
we used \eqref{had},\eqref{fad}
and indicated by $\ts\,\equiv\,$ the equality in
$\Ab$ modulo the ideal $\Jb\ts$.
\qed
\end{proof}

Proposition \ref{prop2} allows us to define for any $c=1\lcd m-1$ 
a linear map
\begin{equation}
\label{xic}
\xic_{\ts c}:\,\Jb\,\ts\backslash\,\Ab\to\Jb\,\ts\backslash\,\Ab
\end{equation}
as the composition $\bar\xi_{\ts c}\,\si_c$ 
applied to the elements of $\Ab$ which were
taken modulo $\Jb\ts$. 

\medskip\smallskip
\noindent\textit{Remark.}
Observe that $\U(\h)\ts\subset\ts{\rm Norm}\ts(\J)\ts$.
Let us denote by $\,\overline{\!{\rm Norm}\ts(\J)\!\!\!}\,\,\,$
the ring of fractions of ${\rm Norm}\ts(\J)$
relative to the same set of denominators \eqref{denset} as before.
Evidently, then $\Jb$ is a two-sided ideal of the ring
$\,\overline{\!{\rm Norm}\ts(\J)\!\!\!}\,\,\,\ts$.
The quotient ring
$$
\Rb\,=\,\Jb\ts\,\backslash\,\,\overline{\!{\rm Norm}\ts(\J)\!\!\!}\,\,\,
$$
bears the same name of \textrm{Mickelsson algebra\/},
as the quotient ring \eqref{malg} does.
One can show \cite{KO} that the linear map \eqref{xic} preserves
the subspace $\Rb\subset\Jb\,\ts\backslash\,\Ab\ts$,
and moreover determines an automorphism of the algebra $\Rb\ts$.
Although we do not use these two facts,
our construction of the $\Y(\gl_n)$-intertwining operator
\eqref{distoper} is underlied by them.
\qed

\medskip\smallskip
In their present form, the operators $\xic_{\ts1}\lcd\xic_{\ts m-1}$  
on the vector space $\Jb\,\ts\backslash\,\Ab$
have been introduced in \cite{KO}.
We will call them \textit{Zhelobenko operators}.
The next proposition states the key property of these
operators; for the proof of this proposition
see \cite[Section 6]{Z}.

\begin{proposition}
The operators\/ $\xic_{\ts1}\lcd\xic_{\ts m-1}$ 
on\/ $\Jb\,\ts\backslash\,\Ab$ satisfy the braid relations
$$
\begin{aligned}
\xic_{\ts c}\,\xic_{\ts c+1}\,\xic_{\ts c}
&\,=\,
\xic_{\ts c+1}\,\xic_{\ts c}\,\xic_{\ts c+1}
\,\quad\quad\textit{for}\ \quad
c<m-1\ts,
\\
\xic_{\ts b}\,\xic_{\ts c}
&\,=\,
\xic_{\ts c}\,\xic_{\ts b}
\hspace{39pt}
\ \quad\textit{for}\ \quad
b<c-1\ts.
\end{aligned}
$$
\end{proposition}

\begin{corollary}
\label{decindep}
For any reduced decomposition $\si=\si_{c_1}\ldots\ts\si_{c_K}$ 
in\/ $\Sym_m$ the composition
$\xic_{\ts c_1}\ldots\,\xic_{\ts c_K}$
of operators on\/ $\Jb\,\ts\backslash\,\Ab$ 
does not depend on the choice of decomposition of\/ $\si\ts$.
\end{corollary}

Recall that $\n'$ denotes the nilpotent subalgebra of $\gl_m$ spanned
by the elements $E_{ab}$ with $a<b\ts$.
Denote by $\Jp$ the left ideal of the algebra $\Ar$ generated by the
elements of the subalgebra $\np\subset\gl_m\ts$. 
Put $\Jpb=\,\overline{\!\U(\h)\!\!\!}\,\,\,\,\Jp\ts$.
Then $\Jpb$ is a left ideal of $\Ab\ts$.
Consider the image $\si_c(\ts\Jpb\ts)\ts$,
this is again a left ideal of $\Ab$.

\begin{proposition}
\label{prop3}
We have\/ $\bar\xi_{\ts c}(\,\si_c(\Jpb\ts))\ts\subset\ts\Jpb\ns+\ts\Jb\ts$.
\end{proposition}

\begin{proof}
Let $\np_{\ts c}$ be the vector subspace in $\gl_m$ spanned
by all elements $E_{ab}$
where $a<b$ but $(a\com b)\neq(c\com c+1)\ts$. Then
$\si_c(\ts\Jpb\ts)\subset\Ab$ is the left ideal generated
by the subspace $\np_{\ts c}\subset\gl_m\ts$ and by 
the element $F_c=E_{c+1,c}\ts.$ Here we use the embedding of
$\gl_m$ into $\Ab\ts$. 
Observe that
the subspace $\np_{\ts c}\subset\gl_m\ts$ is preserved by the
adjoint action of the element $F_c$ on $\gl_m\ts$.
Thus $\xi_{\ts c}\ts(\ts YX)\in\ts\Jpb\ns+\ts\Jb$
for any $X\in\np_c$ and $Y\in\Ar\ts$,
see the definition \eqref{q1}. 

Note that
$\xi_{\ts c}\ts(\ts YF_c)=\xi_{\ts c}\ts(\ts Y)\ts F_c$
for any $Y\in\Ar\ts$, because 
$\widehat{F}_c(\ts YF_c)=\widehat{F}_c(\ts Y)\ts F_c\ts$.
We shall complete the proof of Proposition~\ref{prop3}
by showing that here
$\xi_{\ts c}\ts(\ts Y)\ts F_c\in\Jb\ts$. Indeed,
\begin{gather*}
\xi_{\ts c}\ts(\ts Y)\ts F_c\,=\,
\sum_{s=0}^\infty\,\,
(\ts s\ts!\,H_c^{\ts(s)}\ts)^{-1}\ts E_c^{\ts s}\,
\widehat{F}_c^{\ts s}(\ts Y)\ts F_c\,=
\\
\sum_{s=0}^\infty\,\,
(\ts s\ts!\,H_c^{\ts(s)}\ts)^{-1}\ts E_c^{\ts s}\,F_c\,
\widehat{F}_c^{\ts s}(\ts Y)
\,-\,
\sum_{s=0}^\infty\,\,
(\ts s\ts!\,H_c^{\ts(s)}\ts)^{-1}\ts E_c^{\ts s}\,
\widehat{F}_c^{\ts s+1}(\ts Y)
\,\equiv
\end{gather*}
\begin{gather*}
\sum_{s=1}^\infty\,\,
(\ts s\ts!\,H_c^{\ts(s)}\ts)^{-1}s\ts(H_c-s+1)\ts E_c^{\ts s-1}
\widehat{F}_c^{\ts s}(\ts Y)
\,-\,
\sum_{s=0}^\infty\,\,
(\ts s\ts!\,H_c^{\ts(s)}\ts)^{-1}\ts E_c^{\ts s}\,
\widehat{F}_c^{\ts s+1}(\ts Y)
\,=
\\
\sum_{s=1}^\infty\,\,
(\ts (s-1)\ts!\,H_c^{\ts(s-1)}\ts)^{-1}\ts E_c^{\ts s-1}
\widehat{F}_c^{\ts s}(\ts Y)
\,-\,
\sum_{s=0}^\infty\,\,
(\ts s\ts!\,H_c^{\ts(s)}\ts)^{-1}\ts E_c^{\ts s}\,
\widehat{F}_c^{\ts s+1}(\ts Y)
\,=\,0\ts.
\end{gather*}
Here we used \eqref{fad}, 
and indicated by $\,\equiv\,$ an equality in
$\Ab$ modulo the right ideal $\Jb\ts$.
\qed
\end{proof}

\smallskip\medskip
Proposition \ref{prop3} implies that for each $c=1\lcd m-1$
the Zhelobenko operator \eqref{xic} induces a linear map
$$
\Jb\,\ts\backslash\,\Ab\,/\,\Jpb\,\to\,\Jb\,\ts\backslash\,\Ab\,/\,\Jpb\,.
$$
Now take a weight $\mu\in\h^\ast$ satisfying \eqref{muass}. 
We shall keep the assumption \eqref{muass} on $\mu$
till the end of this section.
Let $\I_{\ts\mu}$ be the left ideal of the algebra $\Ar$
generated by the elements 
$$
E_{ab}
\quad\text{with}\quad
a<b\ts,
\quad
E_{aa}-\mu_a\ts
\,\quad\text{and}\,\quad
\d_{\ts bk}
$$
for all possible $a,b$ and $k\ts$.
Under the isomorphism of $\Ar$ 
with the tensor product \eqref{tenprod}, 
the ideal $\I_{\ts\mu}$ of $\Ar$ corresponds to the ideal of 
the algebra \eqref{tenprod} generated by the elements 
$$
E_{ab}\ot1
\quad\text{with}\quad
a<b\ts,
\quad
E_{aa}\ot1-\mu_a\ts
\,\quad\text{and}\,\quad
1\ot\d_{\ts bk}
$$
for all possible $a,b$ and $k\ts$.
Indeed, for any $a,b=1\lcd m$
the image of the element $E_{ab}\in\Ar$ in the algebra \eqref{tenprod}
is the sum \eqref{eabact}, which equals $E_{ab}\ot1$ plus elements
divisible on the right by the tensor products of the form
$1\ot\d_{\ts bk}\ts$. But the quotient space of \eqref{tenprod}
with respect to the latter ideal can be naturally identified
with the tensor product
$M_{\ts\mu}\ot\P\ts(\ts\CC^{\ts m}\ot\CC^{\ts n})\ts$.
Using the isomorphism of algebras $\Ar$ and \eqref{tenprod},
%%% Under the action of $\Ar$ via the left multiplication,
the quotient space $\Ar\,/\,\I_{\ts\mu}$ can be then also identified with
$M_{\ts\mu}\ot\P\ts(\ts\CC^{\ts m}\ot\CC^{\ts n})\ts$.

Note that
$\mu(H_c)\notin\ZZ$ for any index $c=1\lcd m-1$ due to \eqref{muass}.
Hence we can define the subspace
\,$\Ib_{\ts\mu}=\,\overline{\!\U(\h)\!\!\!}\,\,\,\,\I_{\ts\mu}$ of $\Ab\ts$.
This subspace is also a left ideal of the algebra $\Ab\ts$.
The quotient space $\Ab\,/\,\Ib_{\ts\mu}$ can be still identified 
with the tensor product
$M_{\ts\mu}\ot\P\ts(\ts\CC^{\ts m}\ot\CC^{\ts n})\ts$.
The quotient of $\Ab$ by $\Ib_{\ts\mu}$ and $\Jb$ 
%%% $\Jb\,\ts\backslash\,\Ab\,/\,\Ib_{\ts\mu}$ 
can be then identified with the space of $\n$-coinvariants,
\begin{equation}
\label{dqci}
\Jb\,\ts\backslash\,\Ab\,/\,\Ib_{\ts\mu}
\,=
(\ts M_{\ts\mu}\ot\P\ts(\ts\CC^{\ts m}\ot\CC^{\ts n}))_{\ts\n}\ts.
\end{equation}

Consider the left ideal of the algebra $\Ar$ generated by all the elements
$\d_{\ts bk}$ where $b=1\lcd m$ and $k=1\lcd n\ts$. 
By the definition \eqref{q1}, the image of this ideal under the
map $\xi_{\ts c}$ is contained in the left ideal of $\Ab$
generated by the same elements. The latter ideal is preserved
by the action on $\Ab$ of the element $\si_c\in\Sym_m\ts$. Note that
by \eqref{saction},
$$
\si_c(\ts\mu+\ep_c)=\si_c\circ\mu\ts.
$$
The property \eqref{q12} and Proposition \ref{prop3} now imply that 
the Zhelobenko operator \eqref{xic} induces a linear map
$$
\Jb\,\ts\backslash\,\Ab\,/\,\Ib_{\ts\mu}
\,\to\,
\Jb\,\ts\backslash\,\Ab\,/\,\Ib_{\,\ts\si_c\ts\circ\ts\mu}\,.
$$
Using the identifications 
\eqref{dqci}, 
the Zhelobenko operator \eqref{xic} induces a linear map
\begin{equation}
\label{indmap}
(\ts M_{\ts\mu}\ot\P\ts(\ts\CC^{\ts m}\ot\CC^{\ts n}))_{\ts\n}
\,\to\,
(\ts M_{\ts\si_c\ts\circ\ts\mu}\ot\P\ts(\ts\CC^{\ts m}\ot\CC^{\ts n}))_{\ts\n}
\,.
\end{equation}

\begin{proposition}
\label{inserted}
For any $s=0,1,2,\ts\ldots$
the map \eqref{indmap} commutes with the action of generator
$T_{ij}^{(s+1)}$ of\/ $\Y(\gl_n)$ on the source
and target vector spaces as the element 
\eqref{combact}. 
\end{proposition}

\begin{proof}
Let $Y$ be the element of the algebra $\Ar$
corresponding to the element $\eqref{combact}$ of
the algebra \eqref{tenprod} under the isomorphism of these
two algebras. The element $Y$ then belongs to the 
centralizer of the subalgebra $\U(\gl_n)$ in $\Ar\ts$.
So the left multiplication in $\Ab$ by $Y$
preserves the right ideal \hbox{\ts$\Jb\subset\Ab\ts$,}
and commutes with the linear map 
$\bar\xi_{\ts c}:\Ab\to\Jb\,\ts\backslash\,\Ab\,$; 
see the definition \eqref{q1}.
This left multiplication also commutes with the action of the
element $\si_c\in\Sym_m$ on $\Ab\ts$, because 
the element $Y$ of the algebra $\Ar$ is $\Sym_m$-invariant.
\qed
\end{proof}

The property \eqref{q11} implies that
the restriction of the linear map \eqref{indmap} to the subspace
of vectors of weight $\la$ is a map
$$
(\ts M_{\ts\mu}\ot\P\ts(\ts\CC^{\ts m}\ot\CC^{\ts n}))_{\ts\n}^{\ts\la}
\,\,\to\,
(\ts M_{\ts\si_c\ts\circ\ts\mu}\ot\P\ts(\ts\CC^{\ts m}
\ot\CC^{\ts n}))_{\ts\n}^{\,\si_c\,\circ\ts\la}\,.
$$
Denote the latter map by $I_{\ts c}\,$,
it commutes with the action of $\Y(\gl_n)\ts$
by Proposition \ref{inserted}.

By choosing a reduced decomposition $\si=\si_{c_1}\ldots\ts\si_{c_K}$ 
and taking the composition of operators $I_{c_1}\ldots\ts I_{c_K}$
we obtain an $\Y(\gl_n)$-intertwining operator
$$
I_{\ts\si}:\,
(\ts M_{\ts\mu}\ot\P\ts(\ts\CC^{\ts m}\ot\CC^{\ts n}))_{\ts\n}^{\ts\la}
\,\,\to\,
(\ts M_{\ts\si\,\circ\ts\mu}\ot\P\ts(\ts\CC^{\ts m}
\ot\CC^{\ts n}))_{\ts\n}^{\,\si\,\circ\ts\la}\,.
$$
It does not depend on the choice of the decomposition of
$\si\in\Sym_m$ due to Corollary \ref{decindep}.
This is the operator \eqref{distoper} 
which we intended to exhibit.
Here we identified the symmetric algebra 
$\S\ts(\ts\CC^{\ts m}\ot\CC^{\ts n})$ with the ring
$\P\ts(\ts\CC^{\ts m}\ot\CC^{\ts n})$ in the same way as we
did in Section~1.

{}From now on we will assume that all labels of the weight $\nu=\la-\mu$
are non-negative integers, otherwise both the source and
target modules in \eqref{distoper} are zero. Let
$(\ts\nu_1\lcd\nu_m)$ be the sequence of these labels. Consider the vector
\begin{equation}
\label{munuvec}
1_\mu\ot x_{11}^{\ts\nu_1}\ldots x_{m1}^{\ts\nu_m}
\,\in\ts
M_{\ts\mu}\ot\P\ts(\ts\CC^{\ts m}\ot\CC^{\ts n})
\end{equation}
Note that $x_{11}^{\ts\nu_1}\ldots x_{m1}^{\ts\nu_m}$
is a \textit{highest\/} vector with respect to the 
action of the Lie algebra $\gl_n$ on
$\P\ts(\ts\CC^{\ts m}\ot\CC^{\ts n})\ts$,
any element $E_{ij}\in\gl_n$ with $i<j$
acts on this vector as zero.
With respect to the action of the Lie algebra $\gl_m$
the vector \eqref{munuvec} is of weight $\la\ts$. 
Denote by $v_{\ts\mu}^{\ts\la}$ the image of the vector \eqref{munuvec} in 
$(\ts M_{\ts\mu}\ot\P\ts(\ts\CC^{\ts m}\ot\CC^{\ts n}))_{\ts\n}^{\ts\la}\,$.

\begin{proposition}
\label{isis}
Under condition \eqref{muass}, the vector
$I_{\ts\si}(\ts v_{\ts\mu}^{\ts\la}\ts)$ equals\/
$v_{\ts\si\,\circ\ts\mu}^{\ts\si\,\circ\ts\la}$ multiplied by
\begin{equation}
\label{isim}
\prod\limits_{\substack{1\le a<b\le m \\ \si(a)>\si(b)}}
\ \prod_{r=1}^{\nu_b}\hspace{16pt}
\frac{\,\mu_a-\mu_b-a+b-r}{\,\la_a-\la_b-a+b+r}\ .
\end{equation}
\end{proposition}

\begin{proof}
It suffices to prove the proposition for $\si=\si_c$ with $c=1\lcd m-1\ts$.
Moreover, it then suffices to consider only the case when $m=2$ and hence
$c=1\ts$. Suppose this is the case.
Then under the identification \eqref{dqci} the vector
$$
v_{\ts\mu}^{\ts\la}
\,\in\ts
(\ts M_{\ts\mu}\ot\P\ts(\ts\CC^{\ts 2}\ot\CC^{\ts n}))_{\ts\n}
$$
gets identified with the image in the quotient space
$\Jb\,\ts\backslash\,\Ab\,/\,\Ib_{\ts\mu}$ of the element
$$
x_{11}^{\ts\nu_1}\,x_{21}^{\ts\nu_2}
\,\in\ts\Ar\subset\Ab\ts.
$$
By applying the transposition $\si_1\in\Sym_2$ to this element we get
$
x_{11}^{\ts\nu_2}\,x_{21}^{\ts\nu_1}\,\in\ts\Ar\ts.
$
By applying the operator $\xi_{1}$ to the latter element we get
the sum of elements of $\Ab\ts$,
\begin{equation}
\label{sumel}
\sum_{s=0}^\infty\,\,
(\ts s\ts!\,H_1^{\ts(s)}\ts)^{-1}\ts E_1^{\ts s}\,
\widehat{F}_1^{\ts s}\ts(\ts x_{11}^{\ts\nu_2}\,x_{21}^{\ts\nu_1})\,.
\end{equation}
In particular, here by the definition of the algebra $\Ar$ we have
$$
\widehat{F}_1(\ts x_{11}^{\ts\nu_2}\,x_{21}^{\ts\nu_1})=
[\,F_1\ts\com\ts x_{11}^{\ts\nu_2}\,x_{21}^{\ts\nu_1}\,]=
\sum_{k=1}^n\,\ts
[\,x_{2k}\,\d_{\ts1k}\ts\com\ts x_{11}^{\ts\nu_2}\,x_{21}^{\ts\nu_1}\,]
=\nu_2\cdot x_{11}^{\ts\nu_2-1}\,x_{21}^{\ts\nu_1+1}\,.
$$
Thus the sum \eqref{sumel} equals
\begin{equation}
\label{sumell}
\sum_{s=0}^{\nu_2}\,\,
(\ts s\ts!\,H_1^{\ts(s)}\ts)^{-1}\ts E_1^{\ts s}
\,\cdot\,
\prod_{r=1}^s\,(\ts\nu_2-r+1\ts)
\,\cdot\,
x_{11}^{\ts\nu_2-s}\,x_{21}^{\ts\nu_1+s}\,.
\end{equation}

Since the element $E_1\in\Ar$
is a generator of the left ideal $\I_{\,\ts\si_1\ts\circ\ts\mu}\ts$,
the image of the sum \eqref{sumell} in the quotient vector space
$\Jb\,\ts\backslash\,\Ab\,/\,\Ib_{\,\ts\si_1\ts\circ\ts\mu}$
coincides with that of
\begin{equation}
\label{sume}
\sum_{s=0}^{\nu_2}\,\,
(\ts s\ts!\,H_1^{\ts(s)}\ts)^{-1}\ts
\,\cdot\,
\prod_{r=1}^s\,(\ts\nu_2-r+1\ts)
\,\cdot\,
\widehat{E}_1^{\ts s}\ts
(\ts x_{11}^{\ts\nu_2-s}\,x_{21}^{\ts\nu_1+s}\ts)\,.
\end{equation}
Here $\widehat{E}_1$ is the operator 
of adjoint action on $\PD\ts(\CC^{\ts 2}\ot\CC^{\ts n})\ts$
corresponding to the element $E_1\in\Ar\ts$.
By \eqref{defar} this operator coincides with the
operator of adjoint action of
$$
\ga\ts(E_1)\,=\ts\sum_{k=1}^n\,
x_{1k}\,\d_{\ts2k}\,.
$$
Therefore the sum \eqref{sume} of elements of the algebra $\Ab$ equals
\begin{equation}
\label{sumelll}
\sum_{s=0}^{\nu_2}\,\,
(\ts s\ts!\,H_1^{\ts(s)}\ts)^{-1}\ts
\,\cdot\,
\prod_{r=1}^s\,(\ts\nu_1+r\ts)\,(\ts\nu_2-r+1\ts)
\,\cdot\,
x_{11}^{\ts\nu_2}\,x_{21}^{\ts\nu_1}\,.
\end{equation}

In the sum \eqref{sumelll} the symbol $H_1^{\ts(s)}$ 
stands for the product in the algebra $\Ar\ts$,
$$
\prod_{r=1}^s\,
(\ts H_1-r+1\ts)
\,=\ts
\prod_{r=1}^s\,
(\ts E_{11}-E_{22}-r+1\ts)
\,.
$$
The operator of adjoint action on $\PD\ts(\CC^{\ts 2}\ot\CC^{\ts n})\ts$
corresponding to the element $H_1\in\Ar$
coincides with the operator of adjoint action of the element
$$
\ga\ts(H_1)\,=\ts\sum_{k=1}^n\,
(\,x_{1k}\,\d_{\ts1k}-x_{2k}\,\d_{\ts2k})
$$
The latter operator acts on $x_{11}^{\ts\nu_2}\,x_{21}^{\ts\nu_1}$
as the multiplication by $\nu_2-\nu_1\ts$. 
Since the elements
$$
E_{11}-\mu_2+1
\ \quad\textrm{and}\ \quad
E_{22}-\mu_1-1
$$
are also generators of the left ideal
$\I_{\,\ts\si_1\ts\circ\ts\mu}\ts$,
the image of the sum \eqref{sumelll} in
$\Jb\,\ts\backslash\,\Ab\,/\,\Ib_{\,\ts\si_1\ts\circ\ts\mu}$
coincides with the image of
$$
\sum_{s=0}^{\nu_2}\,\,
\,\prod_{r=1}^s\,\,
\frac
{(\ts\nu_1+r\ts)\,(\ts\nu_2-r+1\ts)}
{\ts r\,(\ts\mu_2-\mu_1+\nu_2-\nu_1-r-1\ts)}
\,\cdot\,
x_{11}^{\ts\nu_2}\,x_{21}^{\ts\nu_1}\,.
$$
The latter image in the quotient vector space
$
\Jb\,\ts\backslash\,\Ab\,/\,\Ib_{\,\ts\si_1\ts\circ\ts\mu}
$
is identified with the vector
$$
v_{\ts\si_1\,\circ\ts\mu}^{\ts\si_1\,\circ\ts\la}
\,\in\ts
(\ts M_{\ts\si_1\ts\circ\ts\mu}\ot\P\ts(\ts\CC^{\ts 2}\ot\CC^{\ts n}))_{\ts\n}
$$
multiplied by the sum
\begin{equation}
\label{eqone}
\sum_{s=0}^{\nu_2}\,\,
\,\prod_{r=1}^s\,\,
\frac
{(\ts\nu_1+r\ts)\,(\ts\nu_2-r+1\ts)}
{\ts r\,(\ts\mu_2-\mu_1+\nu_2-\nu_1-r-1\ts)}
\,\,\,=\,\,
\prod_{r=1}^{\nu_2}\,\,
\frac{\,\mu_1-\mu_2-r+1}{\,\la_1-\la_2+r+1}\ .
\end{equation}

The equality \eqref{eqone} 
%is implicitly contained in \cite[Subsection 5.1.5]{Zhelobenko book} and 
has been also used in \cite[Theorem 8]{TV1}.
Here is its direct proof.
Recall that
$$
\mu_1-\mu_2+\nu_1-\nu_2=\la_1-\la_2\ts.
$$
Let us write $\ts\nu_1=x\ts$, $\ts\nu_2=d\ts$ and $\ts\la_1-\la_2=t\ts$.
Then the equality \eqref{eqone} takes the form
\begin{equation}
\label{eqtwo}
\sum_{s=0}^{d}\,\,
\,\prod_{r=1}^s\,\,
\frac
{(\ts x+r\ts)\,(\ts d-r+1\ts)}
{\ts r\,(\ts-\,t-r-1\ts)}
\,\,\,=\,\,
\prod_{r=1}^{d}\,\,
\frac{\,t-x+d-r+1}{\,t+r+1}\ .
\end{equation}
We shall prove the latter equality for all $x\com t\in\CC$ 
and all nonegative integers $d\ts$.
Define a family of polynomials $\varphi_{d,t}(x)$ of degree $d$ 
in the variable $x$ with coefficients from the field $\CC(t)$
by the following conditions:
\begin{align*}
\label{ccc}
\varphi_{\ts0,t}(x)&=1\ts,
\\
\varphi_{d,t}(x-1)-\varphi_{d,t}(x)
&=d\,(\ts t+2\ts)^{-1}\,\varphi_{d-1,t+1}(x)
\text\quad{for}\quad d\ge1\ts,
\\
\varphi_{d,t}(-1)&=1\ts.
\end{align*}
These conditions uniquely determine the polynomials $\varphi_{d,t}(x)\ts$.
One can easily check that both sides
of the equality \eqref{eqtwo} satisfy these conditions. 
This proves the equality.
\qed
\end{proof}

\begin{comment}

$$
(-1)^s
\,\prod_{r=1}^s
\frac
{(\ts x+r-1\ts)\,(\ts n-r+1\ts)}
{\ts r\,(\ts\la+r+1\ts)}
 -
(-1)^s\,\,
\,\prod_{r=1}^s
\frac
{(\ts x+r\ts)\,(\ts n-r+1\ts)}
{\ts r\,(\ts\la+r+1\ts)}=
$$
$$=(-1)^s\frac{n}{\la+2}
\prod_{r=2}^s\,\frac
{(x+r-1)(\ts n-r+1\ts)}
{\ts r\,(\ts\la+r+1\ts)}(x-(x+s))=$$
$$=(-1)^{s-1}\frac{n}{\la+2}
\prod_{r=1}^{s-1}\,\frac
{\ts (x+r)((n-1)-r+1)\ts}
{\ts r\,(\ts(\la+1)+r+1\ts)}\ ,$$
which is the term of $\varphi_{n-1,\la+1}$     

\end{comment}

Note that 
the $\Y(\gl_n)$-intertwining operator $I_{\ts\si}$ has been defined only
when the weight $\mu$ satisfies the condition \eqref{muass}. 
We also assume that all labels of the weight
$\nu=\la-\mu$ are non-negative integers.
Recall that the sequence of labels $(\ts\rho_1\lcd\rho_m)$
of the weight $\rho$ is $(\ts0\lcd1-m\ts)\ts$.
For any $t\in\CC$ and $N=0,1,2,\ts\ldots$ let us denote by 
$S_{\ts t}^{\ts N}$
the $\Y(\gl_n)$-module obtained from the standard action
of $\U(\gl_n)$ on $\S^N(\CC^{\ts n})\ts$ by pulling back through 
the homomorphism $\pi_n:\Y(\gl_n)\to\U(\gl_n)\ts$, and then
through the automorphism $\tau_{-t}$ of $\Y(\gl_n)\ts$; see
the definitions \eqref{tauz} and \eqref{pin}.
Using Corollary \ref{bimequiv} and the subsequent remarks, 
we can replace the source and target modules in \eqref{distoper}
by equivalent $\Y(\gl_n)$-modules to get an intertwining
operator between two tensor products of $\Y(\gl_n)$-modules,
\begin{equation}
\label{swnop}
S_{\ts\mu_1+\rho_1}^{\ts\nu_1}
\ot\ldots\ot 
S_{\ts\mu_m+\rho_m}^{\ts\nu_m}
\to\,
S_{\ts\widetilde\mu_1+\widetilde\rho_1}^{\ts\widetilde\nu_1}
\ot\ldots\ot 
S_{\ts\widetilde\mu_m+\widetilde\rho_m}^{\ts\widetilde\nu_m}
\end{equation}
where we write
$$
\widetilde\mu_a=\mu_{\si^{-1}(a)}\ts,
\quad 
\widetilde\nu_a=\nu_{\si^{-1}(a)}\ts,
\quad
\widetilde\rho_a=\rho_{\si^{-1}(a)}
$$
for any $a=1\lcd m\ts$.
It is well known that under the condition \eqref{muass} on
the sequence $(\ts\mu_1\lcd\mu_m)$
both tensor products are irreducible $\Y(\gl_n)$-modules,
equivalent to each other\ts; see for instance
\cite[Theorem 3.4]{NT1}. So an intertwining operator 
between these two tensor products is unique up to a multiplier from
$\CC\ts$. For the intertwining operator
corresponding to $I_{\ts\si}$ this multiplier
is determined by Proposition \ref{isis}.
Another expression for an intertwining operator 
between two tensor products of $\Y(\gl_n)$-modules \eqref{swnop}
can be obtained by using a method of I.\,Cherednik \cite{C1},
see for instance \cite[Section 2]{NT1}.

\medskip\smallskip
\noindent\textit{Remark.}
The product \eqref{isim} in Proposition \ref{isis}
does not depend on the choice of reduced decomposition 
$\si_{c_1}\ldots\ts\si_{c_K}$ of the element $\si\in\Sym_m\ts$.
The uniqueness of the
intertwining operator \eqref{swnop} thus provides another proof
of the independence of the composition
$I_{c_1}\ldots\ts I_{c_K}$ on the choice of the
decomposition of $\si\ts$, not involving Corollary \ref{decindep}.
\qed

\medskip\smallskip
A connection between the intertwining operators
on the $m$-fold tensor products of $\Y(\gl_n)$-modules of the 
form $S_{\ts w}^{\ts N}\ts$,
and the results of \cite{Z} for the Lie algebra $\gl_m$
has been established by V.\,Tarasov and A.\,Varchenko \cite{TV2}.
The construction of the operator $I_\si$ given in this section
provides a representation theoretic explanation of that connection. 

%------------------------------------------------------------------------------

\section*{\normalsize 4.\ Olshanski homomorphism}
\setcounter{section}{4}
\setcounter{equation}{0}
\setcounter{theorem}{0}

Let $l$ be a positive integer.
The decomposition $\CC^{\ts n+l}=\CC^{\ts n}\op\CC^{\ts l}$ 
defines an embedding
of the direct sum $\gl_n\op\gl_{\ts l}$ of Lie algebras into $\gl_{n+l}\ts$.
As a subalgebra of $\gl_{n+l}\ts$,
the direct summand $\gl_n$ is spanned by the matrix units 
$E_{ij}\in\gl_{n+l}$ where
$i,j=1\lcd n\ts$. The direct summand $\gl_{\ts l}$ is spanned by
the matrix units $E_{ij}$ where $i,j=n+1\lcd n+l\ts$.
Let $\Cr_{\ts l}$ be the centralizer in 
$\U(\gl_{n+l})$ of the subalgebra $\gl_{\ts l}\subset\gl_{n+l}\ts$.
Set $\Cr_0=\U(\gl_n)\ts$.

Proposition \ref{dast} shows that for any positive integer $m$ 
a homomorphism of associative algebras
\begin{equation}
\label{lastrem}
\Y(\gl_n)\,\to\ts\U(\gl_m)\ot\PD\ts(\CC^{\ts m}\ot\CC^{\ts n})
\end{equation}
can be defined
by mapping $T_{ij}^{\ts(s+1)}$ to the sum \eqref{combact}.
The image of this homomorphism is contained in the
centralizer of the image of $\gl_m$ in 
$\U(\gl_m)\ot\PD\ts(\CC^{\ts m}\ot\CC^{\ts n})\ts$, 
see the remark just after the proof of Proposition \ref{dast}.
In this final section we will compare this homomorphism with
a homomorphism $\Y(\gl_n)\to\Cr_{\ts l}$ defined by G.\,Olshanski \cite{O1}.

Consider the Yangian $\Y(\gl_{n+l})\ts$. 
The subalgebra in $\Y(\gl_{n+l})$ generated by 
$$
T_{ij}^{\ts(1)},T_{ij}^{\ts(2)},\ts\ldots
\quad\text{where}\quad
i\ts,\ns j=1\lcd n 
$$
is isomorphic to $\Y(\gl_n)\ts$ as an associative algebra, see
\cite[Corollary 1.23]{MNO}. Thus we have a natural embedding
$\Y(\gl_n)\to\Y(\gl_{n+l})\ts$, which will be denoted by $\io_{\ts l}\ts$.
Note that $\io_{\ts l}$ is not a Hopf algebra homomorphism.
We also have a surjective homomorphism 
$$
\pi_{n+l}:\ts\Y(\gl_{n+l})\to\U(\gl_{n+l})\ts,
$$
see \eqref{pin}. The composition 
$\pi_{n+l}\,\io_{\ts l}$ coincides with the homomorphism $\pi_n\ts$.

Further, consider the involutive automorphism $\om_{n+l}$
of the algebra $\Y(\gl_{n+l})\ts$, see the definition \eqref{1.51}.
The image of the composition of homomorphisms
$$
\pi_{n+l}\,\ts\om_{n+l}\,\ts\io_{\ts l}:\ts
\Y(\gl_n)\to
%\Y(\gl_{n+l})\to\Y(\gl_{n+l})\to
\U(\gl_{n+l})
$$
belongs to subalgebra $\Cr_{\ts l}\subset\U(\gl_{n+l})\ts$.
Moreover, this image along with the centre
of the algebra $\U(\gl_{n+l})$ generates the subalgebra $\Cr_{\ts l}\ts$.
For the proofs of these two assertions see \cite[Section 2.1]{O2}.
The \textit{Olshanski homomorphism\/} $\Y(\gl_n)\to\Cr_{\ts l}$
is the composition
\begin{equation}
\label{al}
\al_{\ts l}\ts=\ts
\pi_{n+l}\,\ts
\om_{n+l}\,\ts
\io_{\ts l}\,\ts
\tau_{\ts -l}\ts.
\end{equation}

Set $\al_0=\pi_n\ts$.
The family of homomorphisms $\al_0\ts,\al_1,\al_2\ts,\ts\ldots$ has
the following property of stability \cite{O1}.
Let us denote by $\I$ the left ideal of the algebra $\U(\gl_{n+l})$
generated by the elements $E_{1,n+l}\lcd E_{n+l,n+l}\ts$.
One can easily check that
the intersection $\I\cap\Cr_{\ts l}$
is a two-sided ideal of the algebra $\Cr_{\ts l}\ts$.
Moreover, there is a decomposition 
\begin{equation}
\label{aldec}
\Cr_{\ts l}=\Cr_{\ts l-1}\op(\ts\I\cap\Cr_{\ts l})\ts.
\end{equation}
For the proofs of these two assertions again see \cite[Section 2.1]{O2}.
In the equality \eqref{aldec} we regard 
$\Cr_{\ts l-1}\subset\U(\gl_{n+l-1})$ as a subalgebra in 
$\U(\gl_{n+l})$ using the standard embedding
$$
\U(\gl_{n+l-1})\to\U(\gl_{n+l})\ts.
$$
Denote by $\vp_{\ts l}$ the projection to the first direct summand
in \eqref{aldec}. Because $\I\cap\Cr_{\ts l}$
is a two-sided ideal of $\Cr_{\ts l}\ts$, the linear map
$\vp_{\ts l}:\Cr_{\ts l}\to\Cr_{\ts l-1}$ is an algebra homomorphism.

\begin{proposition}
\label{alst}
For any\/ $l=1,2,\ts\ldots$ we have\/ 
$\vp_{\ts l}\,\al_{\ts l}=\al_{\ts l-1}\ts$.
\end{proposition}

\begin{proof}
See \cite[Section 2.1]{O2} once again.
\qed
\end{proof}

Let us describe the homomorphism $\al_{\ts l}$ more explicitly.
For any $i\com j=1\lcd n+l$ regard 
$E_{ij}$ as an element of the algebra $\U(\gl_{n+l})\ts$.
Consider the $(n+l)\times(n+l)$ matrix whose
$ij\ts$-entry is $\de_{ij}-E_{ij}\ts u^{-1}$. 
Consider the inverse of this matrix as a formal power
series in $u^{-1}$ with matrix coefficients.
Denote by $Y_{ij}(u)$ the $ij$-entry of the inverse matrix,
$$
Y_{ij}(u)\ts=\ts
\de_{ij}+E_{ij}\ts u^{-1}+\ts
\sum_{s=1}^\infty\,
\sum_{k_1,\ldots,k_s=1}^{n+l}
E_{ik_1}\ts E_{k_1k_2}\ldots\ts E_{k_{s-1}k_s}\ts E_{\ts k_sj}\,
u^{-s-1}\ts.
$$
Then by \eqref{al},
\begin{equation}
\label{alu}
\al_{\ts l}:\ts T_{ij}(u)\mapsto Y_{ij}(\ts u+l\ts)
\quad\text{for}\quad
i\ts,\ns j=1\lcd n\ts.
\end{equation}
Here each of the formal power series $Y_{ij}(u+l\ts)$ in $(u+l\ts)^{-1}$
should be re-expanded in $u^{-1}$, and %the assignment 
\eqref{alu} is a correspondence
between the respective coefficients of series in $u^{-1}$.

In the present article, we will
employ the homomorphism $\Y(\gl_n)\to\Cr_{\ts l}$
$$
\be_{\ts l}\ts=\ts\al_{\ts l}\,\ts\om_n\ts=\ts
\pi_{n+l}\,\ts
\om_{n+l}\,\ts
\io_{\ts l}\,\ts
\om_n\,\ts
\tau_{\ts l}\ts.
$$
The second equality here follows from the definition \eqref{al}
and from the relation
$$
\tau_{\ts -l}\,\ts\om_n
\ts=\ts
\om_n\,\ts\tau_{\ts l}\ts,
$$
see \eqref{tauz} and \eqref{1.51}. The image
of any series \eqref{tser} under $\be_l$ can be expressed
in terms of quasideterminants \cite[Lemma 4.2]{BK} or
quantum minors \cite[Lemma 8.5]{BK}\ts; see also \cite[Lemma 1.5]{NT2}.
The reason for considering here the homomorphism $\be_{\ts l}$ rather than
$\al_{\ts l}$ will become clear when we state Proposition~\ref{arol}.
Using the definitions of $\be_l$ and $\be_{\ts l-1}$ only,
Proposition \ref{alst} can be restated as

\begin{corollary}
\label{best}
For any\/ $l=1,2,\ts\ldots$ we have\/ 
$\vp_{\ts l}\,\be_{\ts l}=\be_{\ts l-1}\ts$.
\end{corollary}

Consider the Lie algebra $\gl_m$ and its Cartan subalgebra $\h\ts$.
A weight $\mu\in\h^\ast$ is called \textit{polynomial\/} if
its labels $\mu_1\lcd\mu_m$ are non-negative integers such that
$\mu_1\ge\ldots\ge\mu_m\ts$. 
The weight $\mu\in\h^\ast$ is polynomial if and only if
for some non-negative integer $N$ the vector space
$$
\Hom_{\,\gl_m}(\ts L_{\ts\mu}\ts,(\CC^{\ts m})^{\ot N}\ts)\neq\{0\}\ts.
$$
Then 
$$
N=\mu_1+\ldots+\mu_m\ts.
$$
The irreducible $\gl_m$-module $L_{\ts\mu}$ of highest weight $\mu$
is then called a \textit{polynomial\/} module.
Then by setting $\mu_{m+1}=\mu_{m+2}=\ldots=0$
we get a partition 
$(\ts\mu_1,\mu_2\ts,\ts\ldots\ts)$
of $N$. When there is no confusion
with the polynomial weight of $\gl_m\ts$, this partition will be
denoted by $\mu$ as well. The maximal index $a$ with $\mu_a>0$
is then called the \textit{length\/} of the partition, and
is denoted by $\ell(\mu)\ts$. Note that here $\ell(\mu)\le m\ts$.
Further, let $\mup=(\mup_1,\mup_2\ts,\ts\ldots\ts)$ be 
the partition \textit{conjugate\/} to the partition $\mu\ts$.
By definition, here $\mup_b$ is equal to the maximal index $a$ such that
$\mu_a\ge b\ts$. In particular, here $\mup_1=\ell(\mu)\ts$.

Let $\la$ and $\mu$ two polynomial weights of $\gl_m$ such that 
\begin{equation}
\label{lcon}
\ell(\la)\le n+l
\quad\textrm{and}\quad
\ell(\mu)\le l\ts.
\end{equation}
Using the respective partitions, then
$\la$ and $\mu$ can also be regarded as polynomial weights
of the Lie algebras $\gl_{n+l}$ and $\gl_{\ts l}$ respectively.
Denote by $L_{\ts\la}^{\ts\prime}$ and $L_{\ts\mu}^{\ts\prime}$
the corresponding irreducible highest weight modules.
We use this notation to distinguish them from the irreducible modules
$L_{\ts\la}$ and $L_{\ts\mu}$ over the Lie algebra $\gl_m\ts$.

Using the action of the Lie algebra $\gl_{\ts l}$ on 
$L_{\ts\la}^{\ts\prime}$ via its
embedding into $\gl_{n+l}$
as the second direct summand of the subalgebra 
$\gl_n\op\gl_{\ts l}\subset\gl_{n+l}$ consider the vector space
\begin{equation}
\label{hll}
\Hom_{\,\gl_{\ts l}}(\ts L_{\ts\mu}^{\ts\prime}\ts,
L_{\ts\la}^{\ts\prime})\ts.
\end{equation}
The subalgebra $\Cr_{\ts l}\subset\U(\gl_{n+l})$ acts on this vector
space through the action of $\U(\gl_{n+l})$ on $L_{\ts\la}^{\ts\prime}\ts$. 
Moreover, the action of $\Cr_{\ts l}$ on \eqref{hll} is irreducible
\cite[Theorem 9.1.12]{D}. Hence the following identifications
of $\Cr_{\ts l\ts}$-modules are unique up to rescaling of the
vector \text{spaces\ts:}
$$
\Hom_{\,\gl_{\ts l}}(\ts L_{\ts\mu}^{\ts\prime}\ts,
L_{\ts\la}^{\ts\prime})\ts=
\vspace{-4pt}
$$
$$
\Hom_{\,\gl_{\ts l}}(\ts L_{\ts\mu}^{\ts\prime}\ts,
\Hom_{\,\gl_m}(\ts L_{\ts\la}\ts,\S\ts(\ts\CC^{\ts m}\ot\CC^{\ts n+l}\ts)))
\ts=
$$
$$
\Hom_{\,\gl_{\ts l}}(\ts L_{\ts\mu}^{\ts\prime}\ts,
\Hom_{\,\gl_m}(\ts L_{\ts\la}\ts,
\S\ts(\ts\CC^{\ts m}\ot\CC^{\ts l}\ts)\ot
\S\ts(\ts\CC^{\ts m}\ot\CC^{\ts n}\ts)))\ts=
\vspace{2pt}
$$
\begin{equation}
\label{hlls}
\Hom_{\,\gl_m}(\ts L_{\ts\la}\ts,
L_{\ts\mu}\ot\S\ts(\ts\CC^{\ts m}\ot\CC^{\ts n}\ts))\ts.
\vspace{4pt}
\end{equation}
Here we use the classical identifications 
of modules over the Lie algebras $\gl_{n+l}$ and $\gl_m\ts$,
$$
L_{\ts\la}^{\ts\prime}=
\Hom_{\,\gl_m}(\ts L_{\ts\la}\ts,
\S\ts(\ts\CC^{\ts m}\ot\CC^{\ts n+l}\ts))
$$
and
$$
\Hom_{\,\gl_{\ts l}}(\ts L_{\ts\mu}^{\ts\prime}\ts,
\S\ts(\ts\CC^{\ts m}\ot\CC^{\ts l}\ts))=
L_{\ts\mu}
\vspace{4pt}
$$
respectively, see for instance \cite[Section 2.1]{H}.
We also use the decomposition
$$
\S\ts(\ts\CC^{\ts m}\ot\CC^{\ts n+l}\ts)=
\S\ts(\ts\CC^{\ts m}\ot\CC^{\ts l}\ts)\ot
\S\ts(\ts\CC^{\ts m}\ot\CC^{\ts n}\ts)\ts.
$$

By pulling back through the homomorphism $\be_l:\Y(\gl_n)\to\Cr_{\ts l}\ts$,
the vector space \eqref{hlls} becomes a module over the Yangian
$\Y(\gl_n)\ts$. On the other hand, the target vector
space $L_{\ts\mu}\ot\S\ts(\ts\CC^{\ts m}\ot\CC^{\ts n}\ts)$ in \eqref{hlls}
coincides with the vector space of the bimodule $\A_{\ts m}(L_{\ts\mu})$ over
$\gl_m$ and $\Y(\gl_n)$, see the beginning of Section 2. 
Using this bimodule structure, the vector space \eqref{hlls} becomes
another module over $\Y(\gl_n)\ts$.
The next proposition shows that these two
$\Y(\gl_n)$-modules are the same. It also
makes \cite[Remark 12]{A}~more~precise.
We will give a direct proof, 
another proof can be obtained by using \cite[Lemma 4.2]{BK}.

\begin{proposition}
\label{arol}
The action of the algebra\/ $\Y(\gl_n)$ on the vector space \eqref{hlls}
via the homomorphism\/ $\be_l$ coincides with the action inherited from
the bimodule\/ $\A_{\ts m}(L_{\ts\mu})\ts$.
\end{proposition}

\begin{proof}
Consider the action of the subalgebra $\Cr_{\ts l}\subset\U(\gl_{n+l})$
on $\S\ts(\ts\CC^{\ts m}\ot\CC^{\ts n+l}\ts)\ts$.
The Yangian $\Y(\gl_n)$ acts on this vector space via the homomorphism
$\be_l:\Y(\gl_n)\to\Cr_{\ts l}\ts$. Let us identify this vector space with 
$\P\ts(\ts\CC^{\ts m}\ot\CC^{\ts n+l}\ts)$
so that the standard basis vectors of $\CC^{\ts m}\ot\CC^{\ts n+l}$
are identified with the corresponding coordinate functions
$x_{ai}$ 
where
$a=1\lcd m$
and
$i=1\lcd n+l\ts$.
Using the decomposition
\begin{equation}
\label{pmnl}
\P\ts(\ts\CC^{\ts m}\ot\CC^{\ts n+l}\ts)=
\P\ts(\ts\CC^{\ts m}\ot\CC^{\ts l}\ts)\ot
\P\ts(\ts\CC^{\ts m}\ot\CC^{\ts n}\ts)
\vspace{-4pt}
\end{equation}
we will demonstrate that
for any $s=0,1,2,\ts\ldots$ and $i\com j=1\lcd n$ the generator
$T_{ij}^{\ts(s+1)}$ of $\Y(\gl_n)$ then acts on the vector space \eqref{pmnl} 
as the element \eqref{combact} of the tensor product 
$\U(\gl_m)\ot\PD\ts(\CC^{\ts m}\ot\CC^{\ts n})\ts$.
Proposition \ref{arol} will thus follow from Proposition \ref{dast}.

For any $i\com j=1\lcd n+l$
the element $E_{ij}\in\U(\gl_{n+l})$ acts on
$\P\ts(\ts\CC^{\ts m}\ot\CC^{\ts n+l}\ts)$ as the differential operator
$$
\sum_{c=1}^m\,
x_{ci}\,\d_{\ts cj}\,.
$$
Consider the $(n+l)\times(n+l)$ matrix whose $ij$-entry is
$$
\de_{ij}+(u-l)^{-1}\,
\sum_{c=1}^m\,
x_{ci}\,\d_{\ts cj}\,.
$$
Write this matrix and its inverse as the block matrices
$$
\begin{bmatrix}
\,A\,&B\,\\\,C\,&D\,
\end{bmatrix}
\quad\textrm{and}\quad
\begin{bmatrix}\ \At\,&\Bt\ \\ \ \Ct\,&\Dt\ 
\end{bmatrix}
$$
where the blocks $A,B,C,D$ and $\At,\Bt,\Ct,\Dt$ are matrices of sizes
$n\times n$, $n\times l$, $l\times n$, $l\times l$ respectively. 
The action of the algebra $\Y(\gl_n)$ on the vector space
$\P\ts(\ts\CC^{\ts m}\ot\CC^{\ts n+l}\ts)$ via the homomorphism
$\be_l:\Y(\gl_n)\to\Cr_{\ts l}$ can now be described by assigning
to the series $T_{ij}(u)$ with $i\com j=1\lcd n$ the 
$ij$-entry of the matrix $\At^{\,-1}$.

Consider the $(n+l)\times m$ matrix whose $ic\ts$-entry is 
the operator of multiplication by $x_{ci}$ in 
$\P\ts(\ts\CC^{\ts m}\ot\CC^{\ts n+l}\ts)$. Write it as
$$
\begin{bmatrix}
\,P\,
\\
\,\Pb\,
\end{bmatrix}
$$
where the blocks $P$ and $\Pb$ are matrices of sizes
$n\times m$ and $l\times m$ respectively. Similarly,
consider the $m\times(n+l)$ matrix whose $cj$-entry is 
the operator $\d_{cj}\ts$.
Write this matrix as
$$
\begin{bmatrix}
\,Q\,\,\Qb\,\ts
\end{bmatrix}
$$
where the blocks $Q$ and $\Qb$ are matrices of sizes
$m\times n$ and $m\times l$ respectively. Then
$$
\begin{bmatrix}
\,A\,&B\,\\\,C\,&D\,
\end{bmatrix}
\,=\,
1+(u-l)^{-1}
\begin{bmatrix}
\,P\,
\\
\,\Pb\,
\end{bmatrix}
\begin{bmatrix}
\,Q\,\,\Qb\,
\end{bmatrix}
\,=\,
1+(u-l)^{-1}
\begin{bmatrix}\,
\ P\ts Q\,&\,P\ts\Qb\ 
\\ 
\ \Pb\ts Q\,&\,\Pb\ts \Qb\  
\end{bmatrix}
$$ 
where $1$ stands for the $(n+l)\times(n+l)$ identity matrix.
By using Lemma \ref{lemma1}, we get
$$
\At^{\,-1}=\,A-B\ts D^{-1}\ts C
\,=\,
$$
$$
1+(u-l)^{-1}\ts P\ts Q-
(u-l)^{-2}\ts 
P\ts\Qb\,\bigl(\ts1+(u-l)^{-1}\ts\Pb\ts\Qb\,\bigr)^{-1}\ts\Pb\ts Q
\,=\,
$$
\begin{equation}
\label{pulq}
1+P\ts(\ts u-l+\Qb\ts\Pb\,\bigr)^{-1}\ts Q
\vspace{2pt}
\end{equation}
where $1$ now stands for the $n\times n$ identity matrix. 

Consider the $m\times m$ matrix $u-l+\Qb\ts\Pb$ appearing in the line
\eqref{pulq}. Its $ab$-entry is
$$
\de_{ab}\ts(u-l)\ts+\sum_{k=1}^l\,
\d_{\ts a,n+k}\,x_{\ts b,n+k}
\,=\,
\de_{ab}\ts u\ts+\sum_{k=1}^l\,
x_{\ts b,n+k}\,\d_{\ts a,n+k}\,.
$$
Observe that the last displayed sum over $k=1\lcd l$ corresponds 
to the action of the element $E_{\ts ba}\in\U(\gl_m)$ 
on the first tensor factor in the decomposition \eqref{pmnl}. 
Denote by $Z_{ab}(u)$ the $ab$-entry of the matrix inverse to 
$u-l+\Qb\ts\Pb$. The $ij$-entry of the 
matrix \eqref{pulq} can then be written as the sum
$$
\de_{ij}\ts+\sum_{a,b=1}^m x_{ai}\,Z_{ab}(u)\,\d_{\ts bj}
\,=\,
\de_{ij}\ts+\sum_{a,b=1}^m Z_{ab}(u)\,x_{ai}\,\d_{\ts bj}\,.
$$
Using the observation made in the end of Section 1, we now
complete the proof. 
\qed
\end{proof}

Note that the proof of Proposition \ref{arol} remains
valid in the case $l=0$. In this case we assume that $\gl_{\ts
l}=\{0\}\ts$. Further note that the homomorphism 
$\U(\gl_m)\to\PD\ts(\CC^{\ts m}\ot\CC^{\ts l}\ts)$
corresponding to the action of $\gl_m$ 
on the first tensor factor in the decomposition \eqref{pmnl}
is injective if $l\ge m\ts$. Thus, independently of Proposition~\ref{dast},
our proof of Proposition~\ref{arol} shows that 
for any positive integer $m$ a homomorphism of associative algebras
\eqref{lastrem} can be defined
by mapping $T_{ij}^{\ts(s+1)}$ to the sum \eqref{combact}.

Now for any given partitions $\la$ and $\mu$ consider all
integers $l$ satisfying the conditions \eqref{lcon}.
Then consider the corresponding $\Y(\gl_n)$-modules \eqref{hll}
where the integers $l$ vary. By choosing a positive integer
$m$ such that $\ell(\la)\com\ell(\mu)\le m\ts$
we derive from Proposition~\ref{arol} the following known fact,
cf.\ \cite[Theorem~1.6]{N2}.

\begin{corollary}
For all integers\/ $l$ obeying the conditions 
\eqref{lcon} the\/ $\Y(\gl_n)$-modules~\eqref{hll} are equivalent.
\end{corollary}

Further, for any given polynomial weights $\la$ and $\mu$ of
$\gl_m$ we can choose an integer $l$ large
enough to satisfy the conditions \eqref{lcon}. Then
the algebra $\Cr_{\ts l}$ acts on the vector
space \eqref{hlls} irreducibly, while the central elements of
$\U(\gl_{n+l})$ act on \eqref{hlls} via multiplication by scalars.
Hence Proposition \ref{arol} implies another known fact; cf.\
\cite[Theorem 10]{A}.

\begin{corollary}
The action of the algebra\/ $\Y(\gl_n)$ on the vector space \eqref{hlls}
inherited from the bimodule\/ $\A_{\ts m}(L_{\ts\mu})$ is irreducible
for any polynomial weights $\la$ and $\mu$ of\/ $\gl_m\ts$. 
\end{corollary}

Furthermore, it is well known that the vector space \eqref{hll}
is not zero if and only if 
\begin{equation}
\label{llcon}
\la_a\ge\mu_a
\quad\textrm{and}\quad
\lap_a-\mup_a\le n
\quad\textrm{for every}\quad
a=1\lcd m\ts.
\end{equation}
%see for instance \cite[Section I.5]{Macdonald book}\ts. 
By choosing, for given polynomial weights $\la$ and $\mu$ of $\gl_m\ts$,
an integer $l\/$ satisfying %the conditions 
\eqref{lcon}, and then identifying the vector spaces \eqref{hll} 
and \eqref{hlls}, we get another well known fact.

\begin{corollary}
\label{llmnz}
The space \eqref{hlls}
is not zero if and only if the inequalities \eqref{llcon} hold.
\end{corollary}

For further details on the irreducible representations of the Yangian 
$\Y(\gl_n)$ of the form \eqref{hll} see for instance \cite[Section 4]{M} and
\cite[Section 2]{NT2}.

%------------------------------------------------------------------------------

\end{document}